\documentclass[a4paper,12pt]{article}	
\usepackage{amsmath,amsfonts,amsthm,amstext,amssymb,amscd,amsbsy,bezier}
\usepackage[pctex32]{graphics} %figuras coloridas (pctex32)
\usepackage{epsfig}
\usepackage[all]{xy}

\usepackage{epsfig,afterpage}
\usepackage{graphicx}
\usepackage{amsmath}

\newcommand{\seta}{\rightarrow}

\newtheorem{theorem}{Theorem}[section]
\newtheorem{corollary}{Corollary}[section]
\newtheorem{lemma}{Lemma}[section]
\newtheorem{proposition}{Proposition}[section]
\begin{document}

\bibliographystyle{amsplain}

\title{Periodic points on T - fiber bundles over the circle}
\author{Weslem L. Silva and Rafael M. Souza}

\maketitle

\section{Introduction}

Let $f: M \to M$ be a map and $x \in M$, where $M$ a compact manifold. The point $x$ is called a periodic point of $f$ if there exists $n \in \mathbb{N}$ such that $f^{n}(x) = x,$ in this case $x$ a periodic point of $f$ of period $n.$  The set of all 
$\{x \in M| \hbox{ x is periodic}  \}$ is called the set of periodic points of $f$ and is denoted by $P(f).$ 

If $M$ is a compact manifold then the Nielsen theory can be generalized to periodic points. Boju Jiang introduced (Chapter 3 in \cite{Jiang} ) a Nielsen-type homotopy invariant $NF_{n}(f)$ being a lower bound of the number of n-periodic points, for each $g$ homotopic to $f;$
$ Fix(g^{n})  \geq NF_{n}(f).$ 
In case $dim(M) \geq 3$, $M$ compact PL- manifold, then any map $f: M \to M$ is homotopic to a map $g$ satisfying 
$ Fix(g^{n})  = NF_{n}(f) $, this was proved in \cite{Jezierski}. 

Consider a fiber bundle $F \to M \stackrel{p}{\to} B$ where $F,M, B$ are closed manifolds and 
$f: M \to M$ a fiber-preserving map over $B.$ In natural way is to study periodic points of $f$ on $M,$ that is, 
given $n \in \mathbb{N}$ we want to study the set $\{x \in M| f^{n}(x) = x  \}.$ The our main question is; 
when $f$ can be deformed by a fiberwise homotopy to a map  $g: M \to M$ such that $Fix(g^{n}) = \emptyset$  ?
 
In this paper we study periodic points of $f$ when the fiber is the torus $T$ and the base is the circle $S^{1}.$

\section{General problem} 

Let $F \to M \stackrel{p}{\to} B$ be a fibration and $f: M \to M$ a fiber-preserving map over $B$, where 
$F,M,B$ are closed manifolds. Given $n \in \mathbb{N}$, from relation $p \circ f = p,$  we obtain $p \circ f^{n} = p$, thus 
$f^{n}: M \to M$ is also a fiber-preserving map for each $n \in \mathbb{N}$. We want to know 
when $f$ can be deformed by a fiberwise homotopy to a map  $g: M \to M$ such that $Fix(g^{n}) = \emptyset.$
The the following lemma give us a necessary condition to a positive answer the question above.  

\begin{lemma} \label{lemma1}
Let $f: M \to M$ be a fiber-preserving map. If some $k$, where $k$ divides $n$, the map $f^{k}: M \to M$ can not be deformed to a fixed point free map, by a fiberwise homotopy, then can not exists, $g \sim_{B} f ,$ such that $g^{n}: M \to M$ is a fixed point free map.  
\end{lemma}
\proof In fact, suppose that exists $g \sim_{B} f$ such that  $Fix(g^{n}) = \emptyset.$ 
Since $Fix(g^{k}) \subset Fix(g^{n})$ and $Fix(g^{k}) \neq \emptyset$ then we have a contradiction. \qed
\bigskip

Therefore, a necessary condition is that for all $k$, where $k$ divides $n$, the map $f^{k}: M \to M$ must be deformed by a fiberwise 
homotopy to a fixed point free map over $B.$

\bigskip

Note that for each $n$ the square of the following diagram is commutative;
$$\xymatrix{
\ldots \ar[r] & \pi_{1}(F,x_{0}) \ar[r]^{i_{\#}} \ar[d]^{(f^{n}|_{F})_{\#}} & \pi_{1}(M,x_{0}) \ar[r]^{p_{\#}} \ar[d]^{f^{n}_{\#}} 
& \pi_{1}(B,p(x_{0})) \ar[r] \ar[d]^{Id} & 0 \\
\ldots \ar[r] & \pi_{1}(F,f^{n}(x_{0})) \ar[r]^{i_{\#}}  & \pi_{1}(M,f^{n}(x_{0})) \ar[r]^{p_{\#}} 
& \pi_{1}(B,p(x_{0})) \ar[r]  & 0 \\
}$$
 
Since in our case all spaces are path-connected then we will represent the generators of the groups $\pi_{1}(M,f^{n}(x_{0}))$ for each $n,$ 
with the same letters. The same thing we will do with $\pi_{1}(T,f^{n}(0)).$  

\bigskip

Let ${M \times }_{B}  M$ be the pullback of
$p:M \to B$ by $p:M \to B$ and $p_{i}:{M \times }_{B}  M \to M, i=1,2,\,$ the projections to
the first and the second coordinates, respectively.

The inclusion $M \times_{B} M - \Delta
\hookrightarrow M \times_{B} M,$ where $\Delta $ is
the diagonal in $M \times_{B} M,$ is replaced by the
fiber bundle $q:E_{B}(M) \seta M \times_{B}
M,$ whose fiber is denoted by  $\mathcal{ F}.$ We have
$\pi_{m}(E_{B}(M)) \approx \pi_{m}(M \times_{B}
M -\Delta )$ where $E_{B}(M) =\{(x,\omega) \in B
\times A^{I} | i(x) = \omega(0) \},  $ with $A = M
\times_{B} M,$  $B = M \times_{B} M - \Delta$ and $q$ is given by $q(x,\omega)= \omega(1).$

E. Fadell and S. Husseini in \cite{fadell} studied the problem to deform the map $f^{n},$ for each $n \in \mathbb{N},$ to a fixed point 
free map. They supposed that $dim(F)\geq 3$ and that $F,M,B$ are closed manifolds. The necessary and sufficient condition to deform 
$f^{n}$ is given by the following theorem that the proof can be find in \cite{fadell}.  

\begin{theorem} \label{theorem-fadell}
Given $n \in \mathbb{N}$, the map $f^{n}: M \to M$ is deformable to a fixed point free map if and only if there exists a lift 
$\sigma(n)$ in the following diagram; 

\begin{center}
\begin{equation} \label{diagram1.1}
 \xymatrix {       &   \mathcal{F} \ar[d]    &  \mathcal{F} \ar[d]     \\
                  &   E_{B}(f^{n})  \ar[d]^{q_{f^{n}}} \ar[r]^{\bar{q}_{f^{n}}}      & E_{B}(M) \ar[d]^{q}        \\
         M \ar[r]^{1} \ar@{-->}[ru]^{\sigma(n)}         &     M \ar[r]^-{(1,f^{n})}   &   M \times_{B} M  \\ }
\end{equation}
\end{center}

\noindent where $ E_{B}(f^{n}) \to M$ is the fiber bundle induced from $q$ by $(1,f^{n}).$ 
\end{theorem}

In the Theorem \ref{theorem-fadell} we have $\pi_{j-1}(\mathcal{F}) \cong \pi_{j}( M \times_{B} M,  M \times_{B} M - \Delta) $ 
$ \cong \pi_{j}(F,F-x)$ where $x$ is a point in $F.$ In this situation, that is, $dim(F) \geq 3$ the classical obstruction was 
used to find a cross section. 

When $F$ is a surface with Euler characteristic $ \leq 0$ then by Proposition 1.6 from \cite{daci1} 
we have necessary e sufficient conditions to deform $f^{n}$ to a fixed point free map over $B.$
The next proposition gives a relation between a geometric diagram and our problem. 

\begin{proposition} \label{prop-1}
Let $f: M \to M$ be a fiber-preserving map over $B.$ Then there is a map $g$, $g \sim_{B} f,$ 
such that $Fix(g^{n}) = \emptyset$ if and only if there is a map $h_{n}:M \to M \times_{B} M - \Delta$ 
of the form $h_{n} = (Id,s^{n})$, where $s: M \to M,$ is fiberwise homotopic to $f$ and makes the diagram below commutative up to homotopy.
\begin{center}
\begin{equation}\label{diagrama1.1}
\xymatrix{      &       M \times_{B} M - \Delta  \ar @{_{(}->}[d] ^{i}   \\
            M \ar[r]_-{(1, f^{n})} \ar @{-->} [ru]^{h_{n}}    &       M \times_{B} M  \\ }
\end{equation}
\end{center}
\end{proposition}
{\it Proof.} $(\Rightarrow)$ Suppose that exists $g: M \to M$, $g \sim_{B} f,$ with $Fix(g^{n}) = \emptyset.$ Is enough to define 
$h_{n} = (Id, g^{n}),$ that is, $s = g.$ $\phantom{)}$ 

$(\Leftarrow)$ If there is $h_{n}$ then we have $(Id, f^{n}) \sim (Id, s^{n})$, where $s \sim_{B} f$. Therefore 
$Fix(s^{n}) = \emptyset.$ Thus, takes $s=g.$ \qed

\section{Torus fiber-preserving maps} \label{section-2}

Let $T$ be, the torus, defined as the quotient space $ {\mathbb{R} \times \mathbb{R}}/{\mathbb{Z} \times \mathbb{Z}} $. 
We denote by $(x, y)$ the elements of  $\mathbb{R} \times \mathbb{R}$ and by  
$[(x,y)]$ the elements in T.

Let $MA = \frac{T \times [0,1]}{([(x,y)],0) \sim (\left[A \left(^{x}_{y} \right) \right],1)}  $
be the quotient space, where $A$ is a homeomorphism of $T$ induced by an operator in 
$\mathbb{R}^{2}$ that preserves $ \mathbb{Z} \times \mathbb{Z}$. 
The space $MA$ is a fiber bundle over the circle $S^{1}$ where the fiber is the torus.
For more details on these bundles see \cite{daci1}.

Given a fiber-preserving map $f: MA \to MA $, i.e. $p \circ f = p$ we want to compute the number $Fix(g^{n})$  
for each map $g$ fiberwise homotopic to $f$.

Consider the loops in $MA$  given by;
$a(t) = <[(t,0)],0> $, $b(t) = $ $ <[(0,t)],0> $ and 
$c(t) = <[(0,0)],t> $ for $t \in [0,1]$.
We denote by $B$ the matrix of the homomorphism induced on the fundamental group by the 
restriction of $f$ to the fiber $T.$ From \cite{daci1} we have the following theorem
that provides a relationship between the matrices $A$ and $B$, where 
$$ A = { \left( \begin{array}{cc} a_{1} & a_{3} \\ a_{2} & a_{4} \\ \end{array} \right)}  $$
From \cite{daci1} the induced homomorphism $f_{\#}: \pi_{1}(MA) \to \pi_{1}(MA) $ is given by; 
$f_{\#}(a) = a^{b_{1}} b^{b_{2}} $, $f_{\#}(b) = a^{b_{3}} b^{b_{4}} $, $f_{\#}(c) = a^{c_{1}} b^{c_{2}} c $.  
Thus $$ B = { \left( \begin{array}{cc} b_{1} & b_{3} \\ b_{2} & b_{4} \\ \end{array} \right)} $$

\begin{theorem} \label{daci1theorem}
$(1) \,\, \pi_{1}(MA,0) = \langle a,b,c | [a,b] = 1, cac^{-1} = a^{a_{1}} b^{a_{2}} , 
 cbc^{-1}  $ $ = a^{a_{3}} b^{a_{4}} \rangle $ 

$(2)\, B $ commutes with $A$.

$(3) \, $ If $f$ restricted to the fiber is deformable to a fixed point free map then 
the determinant of $B - I$ is zero, where $I$ is the identity matrix.

$(4) \, $ If $v$ is an eigenvector of $B$ associated to 1 $(for \, B \neq Id)$ then $A(v)$ is
also an eigenvector of $B$ associated to 1.

$(5) \, $ Consider $w = A(v)$ if the pair $v,w$ generators $ \mathbb{Z} \times \mathbb{Z}$, 
otherwise let $w$ be another vector so that $v,w$ span $\mathbb{Z} \times \mathbb{Z}$.
Define the linear operator 
$P:\mathbb{R} \times \mathbb{R} \to \mathbb{R} \times \mathbb{R} $ by $P(v) = \left( ^{1}_{0} \right)$
and  $P(w) = \left( ^{0}_{1} \right)$. Consider an isomorphism of fiber bundles, also denoted by $P$,
$P: MA \to M(A^{1})$ where $A^{1} = P \cdot A \cdot P^{-1}$. Then $MA$ is homeomorphic to $M(A^{1})$ 
over $S^{1}$. Moreover we have one of the cases of the table below with $B^{1} = P \cdot A \cdot P^{-1}$
and $B^{1} \neq Id$, except in case $I$:

{\small \begin{center}
\begin{tabular}{|p{1,4cm}|l|}
\hline
 $ Case \,\, I $ & $ A^{1} = \left( \begin{array}{cc} a_{1} & a_{3} \\ a_{2} & a_{4} \\ \end{array} \right) $ ,
 $ B^{1} = \left( \begin{array}{cc} 1 & 0 \\ 0 & 1 \\ \end{array} \right) $ \\
  & $a_{3} \neq 0 $    \\
\hline
$ Case \,\, II $ & $ A^{1} = \left( \begin{array}{cc} 1 & a_{3}\\ 0 & 1 \\ \end{array} \right) $ ,
 $ B^{1} = \left( \begin{array}{cc} 1 & b_{3} \\ 0 & b_{4} \\ \end{array} \right) $ \\
  & $a_{3}(b_{4}-1) = 0 $    \\
\hline
$ Case \,\, III $ & $ A^{1} = \left( \begin{array}{cc} 1 & a_{3} \\ 0 & -1 \\ \end{array} \right) $ ,
 $ B^{1} = \left( \begin{array}{cc} 1 & b_{3} \\ 0 & b_{4} \\ \end{array} \right) $ \\
  & $a_{3}(b_{4}-1) = -2b_{3} $    \\
\hline
$ Case \,\, IV $ & $ A^{1} = \left( \begin{array}{cc} -1 & a_{3} \\ 0 & -1 \\ \end{array} \right) $ ,
 $ B^{1} = \left( \begin{array}{cc} 1 & b_{3} \\ 0 & b_{4} \\ \end{array} \right) $ \\
  & $a_{3}(b_{4}-1) = 0 $    \\  
\hline
$ Case \,\, V $ & $ A^{1} = \left( \begin{array}{cc} -1 & a_{3} \\ 0 & 1 \\ \end{array} \right) $ ,
 $ B^{1} = \left( \begin{array}{cc} 1 & b_{3} \\ 0 & b_{4} \\ \end{array} \right) $ \\
  & $a_{3}(b_{4}-1) = 2b_{3} $    \\  
\hline
\end{tabular} \end{center}} 
\end{theorem}

\bigskip

From \cite{daci1} we have the following theorem:

\begin{theorem} \label{main-theorem-daci1}
A fiber-preserving map $f : MA \to MA$ can be deformed
to a fixed point free map by a homotopy over $S^{1}$ if and only if one
of the cases below holds:

(1) $MA$ is as in case I and $f$ is arbitrary

(2) $MA$ is as in one of the cases $II$ or $III$ and $c_{1}(b_{4}-1)-c_{2}b_{3}=0$

(3) $MA$ is as in case $IV$ and $n(b_{4}(b_{3}+1) -1 -c_{1}(b_{4}-1) + b_{3}c_{2}) - (n-1)(b_{4}-1) \equiv_{2} 0$ except when:

\noindent $a_{3}$ is odd and $[(c_{1}+b_{3}c_{2}\dfrac{(n-1)}{2},nc_{2})]=[(0,0)]\in\dfrac{\mathbb{Z}\oplus\mathbb{Z}}{\langle(1,2),(0,4)\rangle}$ or

\noindent $a_{3}$ is even and $[(nc_{1} +\dfrac{n(n-1)}{2}b_{3}b_{4}c_{2}, c_{2}+(n-1)b_{4}c_{2})]=[(0,0)]$, with $[(0,0)]\in\dfrac{\mathbb{Z}\oplus\mathbb{Z}}{\langle(2,0),(0,2)\rangle}.$

(4) $MA$ is as in case $V$ and either

$a_{3}$ is even and $(b_{4}-1)(c_{1}-\frac{a_{3}}{2}c_{2}-1) \equiv 0 \,\, mod \,\, 2$, except when 
$c_{1}-\frac{a_{3}}{2}c_{2}-1$ and $\frac{b_{4}-1}{L}$ are odd, or 

$a_{3}$ is odd and $\frac{b_{4}-1}{2}(1+c_{2}) \equiv 0 \,\, mod \,\, 2$ except when 
$1+c_{2}$ and $\frac{b_{4}-1}{L}$ are odd, where $L:= mdc(b_{4}-1,c_{2})$.

\end{theorem}

%%%%%%
\bigskip

Given $n \in \mathbb{N}$ we denote the induced homomorphism $f^{n}_{\#}: \pi_{1}(MA) \to \pi_{1}(MA) $ by; 
$f_{\#}(a) = a^{b_{1n}} b^{b_{2n}} $, $f_{\#}(b) = a^{b_{3n}} b^{b_{4n}} $, $f_{\#}(c) = a^{c_{1n}} b^{c_{2n}} c $,   
where $b_{j1} = b_{j}, j=1,...,4$ and $c_{j1}=c_{j}, j=1,2.$ 
Thus the matrix of the homomorphism induced on the fundamental group by the 
restriction of $f^{n}$ to the fiber $T$ is given by;  
$$ B_{n} = { \left( \begin{array}{cc} b_{1n} & b_{3n} \\ b_{2n} & b_{4n} \\ \end{array} \right)} $$
where $B_{1} = B$ is the matrix of $(f_{|T})_{\#}$ and $B_{n} = B^{n}.$ From \cite{K-94} we have
$$ N(h^{n}) = |L(h^{n})| = |det([h_{\#}]^{n} -I)| $$
for each map $h: T \to T$ on torus, where $[h_{\#}]$ is the matrix of induced homomorphism and $I$ is the identity. 

Since $(B^{n}-I) = (B-I)(B^{n-1}+...+B+I)$ then $det(B^{n}-I) = det(B-I)det(B^{n-1}+...+B+I)$. Therefore if $f_{|T}$ is deformable to a fixed point free map then $f^{n}_{|T}$ is deformable to a fixed point free map.

In the Theorems \ref{daci1theorem} and \ref{main-theorem-daci1} putting $f^{n}$ in the place of $f$ we will get conditions 
to $f^{n}.$ The conditions in Theorem \ref{daci1theorem} to $f^{n}$ is the same of $f$ but the conditions to $f^{n}$ 
in the Theorem \ref{main-theorem-daci1} are different of $f$ and are in the Theorem \ref{theorem-1}.

\section{Fixed points of $f^{n}$}

Given a fiber-preserving map $f: MA \to MA$, if $f \sim_{S^{1}} g$ then $f^{n} \sim_{S^{1}} g^{n}$. Therefore, if $Fix(g^{n}) = \emptyset$ 
then the homomorphism $f^{n}_{\#}: \pi_{1}(M) \to \pi(M)$ satisfies the condition of deformability gives in \cite{daci1}.

\begin{proposition} \label{proposition-1}
Let $f: MA \to MA$ be a fiber-preserving map, where $M$ is a T-bundle over $S^{1}.$ Suppose which $f$ restrict to the fiber can 
be deformed to a fixed point free map. This implies $L(f|_{T})=0.$ From Theorem \ref{daci1theorem} we can suppose that the induced homomorphism 
$f_{\#}: \pi_{1}(MA) \to \pi_{1}(MA) $ is given by; $f_{\#}(a) = a $, $f_{\#}(b) = a^{b_{3}} b^{b_{4}} $, $f_{\#}(c) = a^{c_{1}} b^{c_{2}} c .$ Given $n \in \mathbb{N}$ then from relation $(f_{\#})^{n} = f^{n}_{\#}$  we obtain;
$$f^{n}_{\#}(a)=a,$$ $$f^{n}_{\#}(b)=a^{b_{3}\sum_{i=0}^{n-1}b_{4}^{i}}b^{b_{4}^{n}},$$  
$$f^{n}_{\#}(c)=a^{nc_{1}+b_{3}c_{2}\sum_{i=0}^{n-2}(n-1-i)b_{4}^{i}}b^{c_{2}\sum_{i=0}^{n-1}b_{4}^{i}}c.$$
\end{proposition}
{\it Proof.} In fact, $f^{2}_{\#}(b) = f_{\#}(a^{b_{3}} b^{b_{4}}) = a^{b_{3}}(a^{b_{3}} b^{b_{4}})^{b_{4}} = a^{b_{3}+ b_{3}b_{4}}b^{b_{4}^{2}}$ and $f^{2}_{\#}(c) = f_{\#}(a^{c_{1}} b^{c_{2}} c) = a^{c_{1}}(a^{b_{3}} b^{b_{4}})^{c_{2}}(a^{c_{1}} b^{c_{2}} c) = a^{2c_{1}+b_{3}c_{2}}b^{c_{2} + c_{2}b_{4}}c$. Suppose that $f^{n}_{\#}(b) = a^{b_{3}\sum_{i=0}^{n-1}b_{4}^{i}}b^{b_{4}^{n}}$ and $f^{n}_{\#}(c) = a^{nc_{1}+b_{3}c_{2}\sum_{i=0}^{n-2}(n-1-i)b_{4}^{i}}b^{c_{2}\sum_{i=0}^{n-1}b_{4}^{i}}c$. Then, 
\begin{displaymath}
	\begin{array}{lclcl}
		f^{n+1}_{\#}(b) & = & f_{\#}(a^{b_{3}\sum_{i=0}^{n-1}b_{4}^{i}}b^{b_{4}^{n}})                                          & = & a^{b_{3}\sum_{i=0}^{n-1}b_{4}^{i}}(a^{b_{3}} b^{b_{4}})^{b_{4}^{n}} \\
		                & = & a^{b_{3}\sum_{i=0}^{n-1}b_{4}^{i}}(a^{b_{3}b_{4}^{n}} b^{b_{4}^{n+1}})                           & = & a^{b_{3}\sum_{i=0}^{n}b_{4}^{i}}b^{b_{4}^{n+1}};
	\end{array}
\end{displaymath}
\begin{displaymath}
	\begin{array}{lcl}
		f^{n+1}_{\#}(c) & = & f_{\#}(a^{nc_{1}+b_{3}c_{2}\sum_{i=0}^{n-2}(n-1-i)b_{4}^{i}}b^{c_{2}\sum_{i=0}^{n-1}b_{4}^{i}}c) \\
		                & = & a^{nc_{1}+b_{3}c_{2}\sum_{i=0}^{n-2}(n-1-i)b_{4}^{i}}(a^{b_{3}} b^{b_{4}})^{c_{2}\sum_{i=0}^{n-1}b_{4}^{i}}(a^{c_{1}} b^{c_{2}} c)\\
										& = & a^{(nc_{1}+b_{3}c_{2}\sum_{i=0}^{n-2}(n-1-i)b_{4}^{i})+ (b_{3}c_{2}\sum_{i=0}^{n-1}b_{4}^{i})+ (c_{1})} b^{(c_{2}\sum_{i=1}^{n}b_{4}^{i})+(c_{2})} c\\
										& = & a^{(n+1)c_{1}+b_{3}c_{2}\sum_{i=0}^{n-1}(n-1-i)b_{4}^{i}} b^{c_{2}\sum_{i=0}^{n}b_{4}^{i}} c.	
	\end{array}
\end{displaymath}
\qed

%From proposition above and $n = 2$ we have $b_{1}(n) = 1, $ $b_{2}(n) = 0,$ $b_{3}(n) = b_{3}[1+b_{4}],$ 
%$b_{4}(n) = b_{4}^{2},$ $c_{1}(n) = 2c_{1}+c_{2}b_{3}$ and $c_{2}(n) = c_{2}[1+b_{4}].$ 

\begin{theorem} \label{theorem-1}
Let $f: MA \to MA$ be a fiber-preserving map, where $MA$ is a T-bundle over $S^{1}.$ Suppose which $f$ restrict to the fiber can 
be deformed to a fixed point free map and that the induced homomorphism 
$f_{\#}: \pi_{1}(MA) \to \pi_{1}(MA) $ is given by; $f_{\#}(a) = a $, $f_{\#}(b) = a^{b_{3}} b^{b_{4}} $, $f_{\#}(c) = a^{c_{1}} b^{c_{2}} c $  as in cases of the Theorem \ref{main-theorem-daci1}. 
Then $f^{n}: MA \to MA$ can be deformed to a fixed point free map over $S^{1}$ if and only if the following conditions are satisfies;
\bigskip

1) $MA$ is as in case $I$ and $f$ is arbitrary. \bigskip

2) $MA$ is as in cases $II$, $III$ and  $c_{1}(b_{4}-1)-c_{2}b_{3} = 0$ or  $n$ even and $b_{4} = -1.$ \bigskip

3) $MA$ is as in case $IV$ and $n(b_{4}(b_{3}+1) -1 -c_{1}(b_{4}-1) + b_{3}c_{2}) - (n-1)(b_{4}-1) \equiv_{2} 0$ except when:

\noindent $a_{3}$ is odd and $[(nc_{1}+\frac{n(n-1)}{2}b_{3}c_{2},nc_{2})]=[(0,0)]\in\frac{\mathbb{Z}\oplus\mathbb{Z}}{\langle(1,2),(0,4)\rangle}$ or

\noindent $a_{3}$ is even and $[(nc_{1} +\frac{n(n-1)}{2}b_{3}b_{4}c_{2}, c_{2}+(n-1)b_{4}c_{2})]=[(0,0)]\in\frac{\mathbb{Z}\oplus\mathbb{Z}}{\langle(2,0),(0,2)\rangle}.$

\bigskip

4) $MA$ is as in case $V$ and either

$a_{3}$ is even and $n(b_{4}-1)(c_{1}-\frac{a_{3}}{2}c_{2}-1) + (n-1)(b_{4}-1)\equiv 0 \,\, mod \,\, 2$, except when 
$n(c_{1}-\frac{a_{3}}{2}c_{2}-1) + (n-1)$ and $\frac{b_{4}-1}{L}$ are odd, or

$a_{3}$ is odd and $\frac{b_{4}-1}{2}((1+c_{2})(1+(n-1)b_{4})) \equiv 0 \,\, mod \,\, 2$ except when 
$(1+c_{2})(1+(n-1)b_{4})$ and $\frac{b_{4}-1}{L}$ are odd, where $L:= mdc(b_{4}-1,c_{2})$.

\end{theorem}

\bigskip
{\it Proof.} 
In proof of this theorem we will denote $f^{n}_{\#}(b)=a^{b'_{3}}b^{b'_{4}}$ and $f^{n}_{\#}(c)=a^{c'_{1}}b^{c'_{2}}c.$

(1) From Theorem \ref{main-theorem-daci1} each map $f: MA \to MA$ is fiberwise homotopic to a fixed point free map over $S^{1}$ 
in particular that happens with $f^{n}: MA \to MA$ for each $n \in \mathbb{N}.$ 

%$g(<x,y,t>) = <x+c_{1}t+\epsilon, y+c_{2}t+\delta, t>$. Note that $g^{n}(<x,y,t>) = 
%<x+nc_{1}t+n\epsilon, y+nc_{2}t+n\delta, t>$. If $c_{1}c_{2}=0$ then take $\epsilon \in \mathbb{R}-\mathbb{Q}$ 
%or $\delta \in \mathbb{R}-\mathbb{Q}$. In this case we have $Fix(g^{n}) = \emptyset$. 
%If $c_{1}c_{2}\neq 0$ then take $\epsilon \in \mathbb{R}-\mathbb{Q}$ 
%and $\delta = 0$ to get $Fix(g^{n}) = \emptyset$ for each $n \in \mathbb{N}.$

\bigskip
(2) \, If $b_{4}= 1$ then 
$b'_{3}=nb_{3}$, $b'_{4}=1$, 
$c'_{1}=nc_{1}+b_{3}c_{2}\dfrac{n(n-1)}{2}$ and $c'_{2}=nc_{2}$. In this sense, following Theorem \ref{main-theorem-daci1} of \cite{daci1}, in cases $II$ and $III$, $f^{n}$ can be deformed, by a fiberwise homotopy, to a fixed point free map if and only if $c'_{1}(b'_{4}-1)-c'_{2}b'_{3}=0$. However, $c'_{1}(b'_{4}-1)-c'_{2}b'_{3} \ = \ n^{2}c_{2}b_{3}$, and $n^{2}c_{2}b_{3}=0$ if and only if $c_{2}=0$ or $b_{3}=0$.

For $b_{4}\neq 1$ we have $b'_{3}=b_{3}\displaystyle{\sum_{i=0}^{n-1}}b_{4}^{i} = b_{3}\dfrac{b_{4}^{n}-1}{b_{4}-1}$, 
$b'_{4}=b_{4}^{n}$, 
$c'_{1}=nc_{1}+b_{3}c_{2}\displaystyle{\sum_{i=0}^{n-2}}(n-1-i)b_{4}^{i} = nc_{1}+b_{3}c_{2}\dfrac{b_{4}^{n}-nb_{4}+n-1}{(b_{4}-1)^{2}}$ and
$c'_{2}=c_{2}\displaystyle{\sum_{i=0}^{n-1}}b_{4}^{i} =  c_{2}\dfrac{b_{4}^{n}-1}{b_{4}-1}$. Then,
$ c'_{1}(b'_{4}-1)-c'_{2}b'_{3} \ = \ \dfrac{n(b_{4}^{n}-1).(c_{1}(b_{4}-1)-c_{2}b_{3})}{b_{4}-1}.$
In fact,
\begin{displaymath}
	\begin{array}{rcl}
		c'_{1}(b'_{4}-1) & = & \left( nc_{1}+c_{2}b_{3}\dfrac{b_{4}^{n}-nb_{4}+n-1}{(b_{4}-1)^{2}}\right)	(b_{4}^{n}-1)\\
							       & = & nc_{1}(b_{4}^{n}-1) + c_{2}b_{3}\left(\dfrac{(b_{4}^{n}-1)-n(b_{4}-1)}{(b_{4}-1)^{2}}\right)	(b_{4}^{n}-1)\\
										 & = & nc_{1}(b_{4}^{n}-1) + c_{2}b_{3} \left(\dfrac{b_{4}^{n}-1}{b_{4}-1} \right)^{2} -n c_{2}b_{3}\left(\dfrac{b_{4}^{n}-1}{b_{4}-1} \right);\\
		c'_{2}b'_{3}     & = & \left(c_{2}\dfrac{b_{4}^{n}-1}{b_{4}-1} \right) \left( b_{3}\dfrac{b_{4}^{n}-1}{b_{4}-1}\right)
 =  c_{2}b_{3} \left(\dfrac{b_{4}^{n}-1}{b_{4}-1} \right)^{2}.\\								
	\end{array}
\end{displaymath}
Therefore,
\begin{displaymath}
	\begin{array}{ccl}
c'_{1}(b'_{4}-1)-c'_{2}b'_{3} & = & nc_{1}(b_{4}^{n}-1) -n c_{2}b_{3}\left(\dfrac{b_{4}^{n}-1}{b_{4}-1} \right)\\
                              & = & n (b_{4}^{n}-1)\left(c_{1} - \dfrac{c_{2}b_{3}}{b_{4}-1} \right)\\
                              & = & n \left(\dfrac{b_{4}^{n}-1}{b_{4}-1}\right)(c_{1}(b_{4}-1) - c_{2}b_{3})\\
															& = & n (c_{1}(b_{4}-1) - c_{2}b_{3}) \left(\displaystyle{\sum_{i=0}^{n-1}}b_{4}^{i}\right) .
	\end{array}
\end{displaymath}

Note that $cos\left(\dfrac{2k\pi}{n}\right) +i.sin\left(\dfrac{2k\pi}{n}\right)$, for $k=0,1,\dots,n-1$, are the roots of $b_{4}^{n}-1=0$. So, $+1$ and $-1$ are the only two possible integer solutions $b_{4}^{n}-1=0$. Since $b_{4}^{n}-1\neq 0$ for $n$ odd and $b_{4}\neq1$, we may assume that $\dfrac{n(b_{4}^{n}-1).(c_{1}(b_{4}-1)-c_{2}b_{3})}{b_{4}-1}=0$ if and only if $c_{1}(b_{4}-1)-c_{2}b_{3}=0$. Then, by  Theorem \ref{main-theorem-daci1} again, in cases $II$ and $III$ and $n$ odd, $f^{n}$ can be deformed, by a fiberwise homotopy, to a fixed point free map if and only if $f$ can be deformed, by a fiberwise homotopy, to a fixed point free map.
For $n$ even and $b_{4}\neq 1$, $b_{4}^{n}-1 = 0$ if and only if $b_{4}=-1$. 

\bigskip

(3) Following Theorem \ref{main-theorem-daci1} of \cite{daci1}, in cases $IV$, $f^{n}$ can be deformed, by a fiberwise homotopy, to a fixed point free map iff $b'_{4} (b'_{3} + 1) - 1 - c'_{1} (b'_{4} - 1) + c'_{2} b'_{3} \equiv_{2} 0$ except when $a_{3}$ even and $[(c'_{1},c'_{2})]=[(0,0)] \in\dfrac{\mathbb{Z}\oplus\mathbb{Z}}{\langle(2,0),(0,2)\rangle}$, or $a_{3}$ odd and $[(c'_{1},c'_{2})]=[(0,0)]\in\dfrac{\mathbb{Z}\oplus\mathbb{Z}}{\langle(1,2),(0,4)\rangle}$.	

As in case $(2)$, we have
$-c'_{1}(b'_{4}-1)+c'_{2}b'_{3} = -n(c_{1}(b_{4}-1)-c_{2}b_{3})\left(\displaystyle{\sum_{i=0}^{n-1}b_{4}^{i}}\right)$
and 
$ b'_{4}(b'_{3}+1)-1 = b_{4}^{n}\left(1 + b_{3}\displaystyle{\sum_{i=0}^{n-1}b_{4}^{i}}\right)-1 = (b_{4}^{n}-1) + b_{4}^{n}b_{3}\left(\displaystyle{\sum_{i=0}^{n-1}b_{4}^{i}}\right) = (b_{4}-1)\left(\displaystyle{\sum_{i=0}^{n-1}b_{4}^{i}}\right) + b_{4}^{n}b_{3}\left(\displaystyle{\sum_{i=0}^{n-1}b_{4}^{i}}\right)$. 
Then
\begin{displaymath}
\begin{array}{c}
  b'_{4}(b'_{3}+1)-1-c'_{1}(b'_{4}-1)+c'_{2}b'_{3} \ \  = \\
  (b_{4}-1)\left(\displaystyle{\sum_{i=0}^{n-1}b_{4}^{i}}\right) + b_{3}b_{4}^{n}\left(\displaystyle{\sum_{i=0}^{n-1}b_{4}^{i}}\right) -n(c_{1}(b_{4}-1)-c_{2}b_{3})\left(\displaystyle{\sum_{i=0}^{n-1}b_{4}^{i}}\right) \ \ = \\
  (b_{4}^{n}-1)(1-nc_{1}) +b_{3}(b_{4}^{n}+nc_{2})\left(\displaystyle{\sum_{i=0}^{n-1}b_{4}^{i}}\right) \ \ \equiv_{2}\\
  (b_{4}-1)(1-nc_{1}) +b_{3}(b_{4}+nc_{2})(1+(n-1)b_{4}) \ = \\
	(b_{4}-1)(1-nc_{1}) +b_{3}(b_{4}+(n-1)b_{4}^{2} +nc_{2} +n^{2}b_{4}c_{2}-nb_{4}c_{2}) \ \equiv_{2} \\
	(b_{4}-1)(1-nc_{1}) +b_{3}(b_{4}+(n-1)b_{4} +nc_{2} +nb_{4}c_{2}-nb_{4}c_{2}) \ = \\
	(b_{4}-1)(1-nc_{1}) +b_{3}(nb_{4} +nc_{2}) \ = \\
 	n(b_{4}(b_{3}+1) -1 -c_{1}(b_{4}-1) + b_{3}c_{2})-(n-1)(b_{4}-1).
	\end{array}
\end{displaymath}

The exceptions holds for $a_{3}$ even and $[(c'_{1},c'_{2})]=[(0,0)] \in\dfrac{\mathbb{Z}\oplus\mathbb{Z}}{\langle(2,0),(0,2)\rangle}$, or $a_{3}$ odd and $[(c'_{1},c'_{2})]=[(0,0)]\in\dfrac{\mathbb{Z}\oplus\mathbb{Z}}{\langle(1,2),(0,4)\rangle}$. In this sense, we have
$$(c'_{1},c'_{2})= \displaystyle{\left(nc_{1}+b_{3}c_{2}\sum_{i=0}^{n-1}(n-1-i)b_{4}^{i}, c_{2}\sum_{i=0}^{n-1}b_{4}^{i}\right)}.$$
If $a_{3}$ is odd then $b_{4}=1$, $c_{2}\displaystyle{\sum_{i=0}^{n-1}1^{i}}= nc_{2}$ and $nc_{1}+b_{3}c_{2}\displaystyle{\sum_{i=0}^{n-1}(n-1-i)1^{i}} = nc_{1}+b_{3}c_{2}\dfrac{n(n-1)}{2}$. If $a_{3}$ is even then $c_{2}\displaystyle{\sum_{i=0}^{n-1}b_{4}^{i}}\equiv_{2} c_{2}(1+(n-1)b_{4})$ and $nc_{1}+b_{3}c_{2}\displaystyle{\sum_{i=0}^{n-1}(n-1-i)b_{4}^{i} }\equiv_{2}nc_{1}+ \dfrac{n(n-1)}{2}b_{3}b_{4}c_{2}$.

\bigskip

(4) From Theorem \ref{main-theorem-daci1} the map $f^{n}$ can be deformed, over $S^{1}$, to a fixed point free map if and only if 
the following condition is satisfy;

$a_{3}$ is even and $(b^{'}_{4}-1)(c^{'}_{1}-\frac{a_{3}}{2}c^{'}_{2}-1) \equiv 0 \,\, mod \,\, 2$, except when 
$c^{'}_{1}-\frac{a_{3}}{2}c^{'}_{2}-1$ and $\frac{b^{'}_{4}-1}{L}$ are odd, or 
$a_{3}$ is odd and $\frac{b^{'}_{4}-1}{2}(1+c^{'}_{2}) \equiv 0 \,\, mod \,\, 2$ except when 
$1+c^{'}_{2}$ and $\frac{b^{'}_{4}-1}{L}$ are odd, where $L:= mdc(b^{'}_{4}-1,c^{'}_{2})$.

Note that if $b_{4}=1$ then from Theorem \ref{daci1theorem} we must have $b_{3}=0$ and this situation return in 
the case $I$. Therefore let us suppose $b_{4} \neq 1.$ 

From previous calculation we have;
$b^{'}_{4} = b^{n}_{4} $, $b^{'}_{3} = b_{3}\frac{b^{n}_{4}-1}{b_{4}-1}$, $c^{'}_{2} = c_{2}\frac{b^{n}_{4}-1}{b_{4}-1}$  
 and $c^{'}_{1} = nc_{1} + b_{3}c_{2}\frac{b^{n}_{4}-nb_{4}+n-1}{(b_{4}-1)^{2}}.$ From Theorem \ref{daci1theorem} we have 
$a_{3}(b_{4}-1) = 2b_{3}.$ 

Suppose $a_{3}$ even. Since $c^{'}_{1}(b^{'}_{4}-1)-c^{'}_{2}b^{'}_{3} = \frac{n(b^{n}_{4}-1)(c_{1}(b_{4}-1)-c_{2}b_{3})}{b_{4}-1} $     
Then $(b^{'}_{4}-1)(c^{'}_{1}-\frac{a_{3}}{2}c^{'}_{2}-1) = n(b^{n}_{4}-1)(c_{1}-\frac{a_{3}}{2}c_{2}-1)+(n-1)(b^{n}_{4}-1).$ 
In fact,
$$\begin{array}{lll}
c^{'}_{1}-\frac{a_{3}}{2}c^{'}_{2} & = & nc_{1} + b_{3}c_{2}\frac{b^{n}_{4}-nb_{4}+n-1}{(b_{4}-1)^{2}} 
- \frac{a_{3}}{2}c_{2}\frac{b^{n}_{4}-1}{b_{4}-1} \\
& = & nc_{1} + b_{3}c_{2}\frac{(b^{n}_{4}-1)-n(b_{4}-1)}{(b_{4}-1)^{2}} - b_{3}c_{2}\frac{b^{n}_{4}-1}{(b_{4}-1)^{2}} \\
& = & nc_{1} - \frac{b_{3}c_{2}n}{b_{4}-1} \\
& = &  n(c_{1}-\frac{a_{3}}{2}c_{2}).\\
\end{array}
$$
 
We know that if $L:=mdc(b_{4}-1,c_{2})$ then $kL:= mdc(k(b_{4}-1),kc_{2}).$ Thus, $mdc(b^{'}_{4}-1,c^{'}_{2}) = L^{'} =
\frac{b^{n}_{4}-1}{(b_{4}-1)}L = $ where $L = mdc(b_{4}-1,c_{2})$ because
$b^{'}_{4}-1 = \frac{b^{n}_{4}-1}{(b_{4}-1)} (b_{4}-1)$ and $c^{'}_{2} =c_{2} \frac{b^{n}_{4}-1}{(b_{4}-1)}.$  
Furthermore, $\frac{b^{'}_{4}-1}{L^{'}} = \frac{b^{'}_{4}-1}{L} \frac{b_{4}-1}{(b^{n}_{4}-1)} = \frac{b_{4}-1}{L}. $
With these calculations we obtain the conditions statements on the theorem.

In the case $a_{3}$ odd we must have; $\frac{b^{n}_{4}-1}{2}(1+c_{2}\frac{b^{n}_{4}-1}{b_{4}-1}) \equiv 0 \,\, mod \,\, 2$ except when 
$1+c_{2}\frac{b^{n}_{4}-1}{b_{4}-1}$ and $\frac{b_{4}-1}{L}$ are odd, where $L:= mdc(b_{4}-1,c_{2})$. 

Note that $\frac{b^{n}_{4}-1}{b_{4}-1}$ is even if and only if $1+(n-1)b_{4}$ is even, and $b^{n}_{4}-1$ is even 
if and only if $b_{4}-1$ is even, for all $n \in \mathbb{N}.$ With this we obtain the enunciate of the theorem. 
\qed
\bigskip

\begin{corollary} \label{corollary-1}
From Theorem \ref{theorem-1}, if $f: MA \to MA$ is deformed to a fixed point free map over $S^{1}$ and $n \in \mathbb{N}$ is odd then 
the map $f^{n}: MA \to MA$ can be deformed to a fixed point free map over $S^{1}.$
\end{corollary}
{\it Proof.} 
If $f: MA \to MA$ is deformed to a fixed point free map over $S^{1}$ then  the conditions of the Theorem \ref{main-theorem-daci1} are satisfied. Suppose $n$ odd then the conditions of the Theorem \ref{theorem-1} also are satisfied. Thus $f^{n}: MA \to MA$ can be deformed to a fixed point free map over $S^{1}.$ \qed

\bigskip

In the corollary above if $n$ is even the above statement may not holds, for example in the case V of the Theorem 
\ref{theorem-1} if $n$, $b_{4}$, $a_{3}$ and $c_{1}- \frac{a_{3}}{2}c_{2} -1$ are even then $f: MA \to MA$ is deformed to a fixed point free map over $S^{1}$ but $f^{n}$ is not.  
 
\begin{proposition} \label{proposition-2}
Let $f: MA \to MA$ be a fiber-preserving such that the induced homomorphism 
$f_{\#}: \pi_{1}(MA) \to \pi_{1}(MA) $ is given by; $f_{\#}(a) = a $, $f_{\#}(b) = a^{b_{3}} b^{b_{4}} $, 
$f_{\#}(c) = a^{c_{1}} b^{c_{2}} c .$ Suppose that for some $n$ odd, $ n \in \mathbb{N}$ the fiber-preserving map 
$f^{n}: MA \to MA$ is deformed to a fixed point free map over $S^{1}.$ If $k$ is a divisor of $n$ then the map $f^{k}: MA \to MA$ can be deformed by a fiberwise homotopy to a fixed point free map over $S^{1}.$
\end{proposition}
{\it Proof.}
Is enough to verify that if the conditions of the Theorem \ref{theorem-1} are satisfied for some $n$ odd then those conditions are also 
satisfied for any $k$ divisor of $n.$ We will analyze each case of the Theorem \ref{theorem-1}. 

\bigskip
Case I. In this case for each $n \in \mathbb{N}$ the fiber-preserving map can be deformed over $S^{1}$ to a fixed point free map. 

\bigskip
Cases II and III. In these cases if for some $n$ odd the fiber-preserving map $f^{n}: MA \to MA$ is deformed to a 
fixed point free map over $S^{1}$ then we must have; $c_{1}(b_{4}-1)-c_{2}b_{3} = 0.$ Thus for all $k \leq n$ the 
$f^{k}$ can be deformed to a fixed point free map over $S^{1},$ in particular when $k$ divides $n.$ 

\bigskip
Case IV. Suppose that for some odd positive integer $n$ the fiber-preserving map $f^{n}: MA \to MA$ is deformed to a fixed point free map over $S^{1}$, then $n(b_{4}(b_{3}+1) -1 -c_{1}(b_{4}-1) + b_{3}c_{2}) - (n-1)(b_{4}-1) \equiv_{2} 0$ and

\noindent if $a_{3}$ is odd then $[(nc_{1}+\frac{n(n-1)}{2}b_{3}c_{2},nc_{2})]\neq[(0,0)]\in\frac{\mathbb{Z}\oplus\mathbb{Z}}{\langle(1,2),(0,4)\rangle}$ or

\noindent if $a_{3}$ is even then $[(nc_{1} +\frac{n(n-1)}{2}b_{3}b_{4}c_{2}, c_{2}+(n-1)b_{4}c_{2})]\neq[(0,0)]\in\frac{\mathbb{Z}\oplus\mathbb{Z}}{\langle(2,0),(0,2)\rangle}.$

Note that if $a_{3}$ is odd then $b_{4}=1$ and
$$n(b_{4}(b_{3}+1) -1 -c_{1}(b_{4}-1) + b_{3}c_{2}) - (n-1)(b_{4}-1) \equiv_{2} b_{3} + b_{3}c_{2}=b_{3}(1 + c_{2}); $$
$$\begin{array}{ccl}
[(nc_{1}+\frac{n(n-1)}{2}b_{3}c_{2},nc_{2})] & = & [(0,nc_{2}- 2(nc_{1}+\frac{n(n-1)}{2}b_{3}c_{2}))]\\
& = & [(0,n(c_{2} - 2c_{1}-(n-1)b_{3}c_{2} )))] \in\frac{\mathbb{Z}\oplus\mathbb{Z}}{\langle(1,2),(0,4)\rangle};\\
&\Rightarrow& n(c_{2} - 2c_{1}-(n-1)b_{3}c_{2})\not\equiv_{4}0\\
&\Rightarrow& c_{2} - 2c_{1}-(n-1)b_{3}c_{2}\not\equiv_{4}0.
\end{array}$$
If $b_{3}$ is even then $c_{2} - 2c_{1}-(n-1)b_{3}c_{2}\equiv_{4}c_{2} - 2c_{1}$. If $c_{2}$ is odd then $c_{2} - 2c_{1}-(n-1)b_{3}c_{2}$ is odd and $c_{2} - 2c_{1}-(k-1)b_{3}c_{2}$ is odd for each odd $k$.  If $a_{3}$ is even we have 
$$\left[\left(nc_{1} +\frac{n(n-1)}{2}b_{3}b_{4}c_{2}, c_{2}+(n-1)b_{4}c_{2}\right)\right]=\left[\left(c_{1} +\frac{(n-1)}{2}b_{3}b_{4}c_{2}, c_{2}\right)\right].$$
Then, $c_{2}\equiv_{2}1$ or $c_{1} +\frac{(n-1)}{2}b_{3}b_{4}c_{2}\equiv_{2}1$. So, if $c_{2}\equiv_{2}0$ then $c_{1} +\frac{(n-1)}{2}b_{3}b_{4}c_{2}\equiv_{2}c_{1}$. Therefore, $c_{1}\equiv_{2}1$ or $c_{2}\equiv_{2}1$.
 
Let $k$ be an integer such that $k$ divides $n$ then $k$ must be odd and
$$\begin{array}{ccc}
k(b_{4}(b_{3}+1) -1 -c_{1}(b_{4}-1) + b_{3}c_{2}) - (k-1)(b_{4}-1) & \equiv_{2} &  \\
n(b_{4}(b_{3}+1) -1 -c_{1}(b_{4}-1) + b_{3}c_{2}) - (n-1)(b_{4}-1) & \equiv_{2} & 0.
\end{array}$$
If $a_{3}$ is odd then $b_{3}$ is even or $c_{2}$ is odd, $b_{4}=1$ and
$$\begin{array}{ccl}
[(kc_{1}+\frac{k(k-1)}{2}b_{3}c_{2},kc_{2})] & = & [(0,k(c_{2} - 2c_{1}-(k-1)b_{3}c_{2} )))] \\
& = & [(0,c_{2} - 2c_{1}-(k-1)b_{3}c_{2} )))] \\
& = & [(0,c_{2} - 2c_{1}-(n-1)b_{3}c_{2} )))] \\
&\neq & [(0,0)] \in\frac{\mathbb{Z}\oplus\mathbb{Z}}{\langle(1,2),(0,4)\rangle}.\\
\end{array}$$
Then, $f^{k}: MA \to MA$ can be deformed to a fixed point free map over $S^{1}$.

If $a_{3}$ is even then
$$\begin{array}{ccl}
[(kc_{1} +\frac{k(k-1)}{2}b_{3}b_{4}c_{2}, c_{2}+(k-1)b_{4}c_{2})] & = & [(c_{1} +\frac{(k-1)}{2}b_{3}b_{4}c_{2}, c_{2})] \\
&\neq & [(0,0)]\in\frac{\mathbb{Z}\oplus\mathbb{Z}}{\langle(2,0),(0,2)\rangle}.\\
\end{array}$$
Then, $f^{k}: MA \to MA$ can be deformed to a fixed point free map over $S^{1}$.

\bigskip
Case V.  Suppose that for some $n$ odd, $ n \in \mathbb{N}$ the fiber-preserving map 
$f^{n}: MA \to MA$ is deformed to a fixed point free map over $S^{1}.$ If $k$ divides $n$ then there is $l \in \mathbb{N}$ 
such that $kl=n,$ in particular $l$ must be odd. The conditions to deform $f^{n}$ in this case with $a_{3}$ even are;

$a_{3}$ is even and $n(b_{4}-1)(c_{1}-\frac{a_{3}}{2}c_{2}-1) + (n-1)(b_{4}-1)\equiv 0 \,\, mod \,\, 2$, except when 
$n(c_{1}-\frac{a_{3}}{2}c_{2}-1) + (n-1)$ and $\frac{b_{4}-1}{L}$ are odd, where $L:=mdc(b_{4}-1,c_{2}).$ 

We have; $l[k(b_{4}-1)(c_{1}-\frac{a_{3}}{2}c_{2}-1) + (k-1)(b_{4}-1)] = $ 
$ lk(b_{4}-1)(c_{1}-\frac{a_{3}}{2}c_{2}-1) + (lk-l)(b_{4}-1) =$
$ n(b_{4}-1)(c_{1}-\frac{a_{3}}{2}c_{2}-1) + (n-1)(b_{4}-1) + (1-l)(b_{4}-1).$ 

From hypothesis $n(b_{4}-1)(c_{1}-\frac{a_{3}}{2}c_{2}-1) + (n-1)(b_{4}-1)$ is even. Since $l$ is odd then $(1-l)(b_{4}-1)$ is even. 
Therefore $l[k(b_{4}-1)(c_{1}-\frac{a_{3}}{2}c_{2}-1) + (k-1)(b_{4}-1)]$ is even. Since $l$ is odd then 
we must have $k(b_{4}-1)(c_{1}-\frac{a_{3}}{2}c_{2}-1) + (k-1)(b_{4}-1) \equiv 0 \,\, mod \,\, 2.$

Note that $l[k(c_{1}-\frac{a_{3}}{2}c_{2}-1) + (k-1)] = $
$[n(c_{1}-\frac{a_{3}}{2}c_{2}-1) + (n-1)]+ (1-l).$ Thus $l[k(c_{1}-\frac{a_{3}}{2}c_{2}-1) + (k-1)]$ is odd because 
$[n(c_{1}-\frac{a_{3}}{2}c_{2}-1) + (n-1)]$ is odd by hypothesis and $(1-l)$ is even. Since $l$ is odd then we must have 
$k(c_{1}-\frac{a_{3}}{2}c_{2}-1) + (k-1)$ odd. Therefore the conditions to deform $f^{k}$ to a fixed point free map over $S^{1}$, in the Theorem \ref{theorem-1}, are satisfied. The case $a_{3}$ odd is analogous.
\qed

\bigskip

\begin{proposition} \label{proposition-3}
Let $f: MA \to MA$ be a fiber-preserving. If $m,n$ are odd, $m,n \in \mathbb{N},$ then $f^{m}$ is deformable to a fixed 
point free map over $S^{1}$ if and only if $f^{n}$ is deformable to a fixed point free map over $S^{1}.$ 
\end{proposition}
{\it Proof.} If $m,n$ are odd and $f^{m}$ is deformable to a fixed point free map over $S^{1}$ then by Proposition \ref{proposition-2} 
$f$ is deformable to a fixed point free map over $S^{1}$. From Corollary \ref{corollary-1}  
$f^{n}$ is deformable to a fixed point free map over $S^{1}.$ \qed

\bigskip

We have a analogous result to $n$ even;

\begin{proposition} \label{proposition-4}

Let $f: MA \to MA$ be a fiber-preserving map, where $MA$ is a T-bundle over $S^{1}.$ Suppose which $f$ restrict to the fiber can be deformed to a fixed point free map and that the induced homomorphism $f_{\#}: \pi_{1}(MA) \to \pi_{1}(MA) $ is given by; $f_{\#}(a) = a $, $f_{\#}(b) = a^{b_{3}} b^{b_{4}} $, $f_{\#}(c) = a^{c_{1}} b^{c_{2}} c $  as in cases of the Theorem \ref{main-theorem-daci1}. Given an even positive integer $n$ such that $f^{n}$ is deformable to a fixed point free map over $S^{1}$ then $f^{k}$ is deformable to a fixed point free map over $S^{1}$, for all even positive integer $k$, except when $MA$ is as in case $IV$ and 

\noindent $a_{3}$ is odd and $k\equiv_{4}0$ or

\noindent $a_{3}$ is even, $k\equiv_{4}0$ and $b_{3}b_{4}c_{2}\equiv_{2}1$.
\end{proposition}

\bigskip

{\it Proof.}
Is enough to verify that if the conditions of the Theorem \ref{theorem-1} are satisfied for some $n$ even then those conditions are also 
satisfied by every even $k$. We will analyze each case of the Theorem \ref{theorem-1}. 

\bigskip

Case I. In this case for each $n \in \mathbb{N}$ the fiber-preserving map can be deformed over $S^{1}$ to a fixed point free map. 

\bigskip
Cases II and III. In these cases if for some $n$ odd the fiber-preserving map $f^{n}: MA \to MA$ is deformed to a 
fixed point free map over $S^{1}$ then we must have; $c_{1}(b_{4}-1)-c_{2}b_{3} = 0$ or $b_{4}=-1$. Thus, for all even $k$, $f^{k}$ can be deformed to a fixed point free map over $S^{1}$.

\bigskip
Case IV. If $n$ is an even positive integer and $f^{n}: MA \to MA$ is deformed to a fixed point free map over $S^{1}$, then $n(b_{4}(b_{3}+1) -1 -c_{1}(b_{4}-1) + b_{3}c_{2}) - (n-1)(b_{4}-1) \equiv_{2} 0$ and

\noindent if $a_{3}$ is odd then $[(nc_{1}+\frac{n(n-1)}{2}b_{3}c_{2},nc_{2})]\neq[(0,0)]\in\frac{\mathbb{Z}\oplus\mathbb{Z}}{\langle(1,2),(0,4)\rangle}$ or

\noindent if $a_{3}$ is even then $[(nc_{1} +\frac{n(n-1)}{2}b_{3}b_{4}c_{2}, c_{2}+(n-1)b_{4}c_{2})]\neq[(0,0)]\in\frac{\mathbb{Z}\oplus\mathbb{Z}}{\langle(2,0),(0,2)\rangle}.$

Note that if $a_{3}$ is odd then $b_{4}=1$ and
$$\begin{array}{ccl}
[(nc_{1}+\frac{n(n-1)}{2}b_{3}c_{2},nc_{2})] & = & [(0,nc_{2}- 2(nc_{1}+\frac{n(n-1)}{2}b_{3}c_{2}))]\\
& = & [(0,n(c_{2} - 2c_{1}-(n-1)b_{3}c_{2} )))] \in\frac{\mathbb{Z}\oplus\mathbb{Z}}{\langle(1,2),(0,4)\rangle};\\
&\Rightarrow& n(c_{2} - 2c_{1}-(n-1)b_{3}c_{2})\not\equiv_{4}0\\
&\Rightarrow& c_{2} - 2c_{1}-(n-1)b_{3}c_{2}\equiv_{2}1 \ and \ n\equiv_{4}2.
\end{array}$$
If $a_{3}$ is even we have 
$$\left[\left(nc_{1} +\frac{n(n-1)}{2}b_{3}b_{4}c_{2}, c_{2}+(n-1)b_{4}c_{2}\right)\right]=\left[\left(\frac{n}{2}b_{3}b_{4}c_{2}, c_{2}(1+b_{4})\right)\right].$$

Let $k$ be an even positive integer then $k$ must be odd and
$$\begin{array}{ccc}
k(b_{4}(b_{3}+1) -1 -c_{1}(b_{4}-1) + b_{3}c_{2}) - (k-1)(b_{4}-1) & \equiv_{2} &  0\\
\end{array}$$
Then, $f^{k}: MA \to MA$ can be deformed to a fixed point free map over $S^{1}$ except when:

\noindent $a_{3}$ is odd and 
$$\begin{array}{ccl}
[(kc_{1}+\frac{k(k-1)}{2}b_{3}c_{2},kc_{2})] & = & [(0,k(c_{2} - 2c_{1}-(k-1)b_{3}c_{2} )))] \\
& = & [(0,k )))] \\
& = & [(0,0)] \in\frac{\mathbb{Z}\oplus\mathbb{Z}}{\langle(1,2),(0,4)\rangle}.\\
\end{array}$$

\noindent or $a_{3}$ is even and
$$\begin{array}{ccl}
[(kc_{1} +\frac{k(k-1)}{2}b_{3}b_{4}c_{2}, c_{2}+(k-1)b_{4}c_{2})] & = & \left[\left(\frac{k}{2}b_{3}b_{4}c_{2}, c_{2}(1+b_{4})\right)\right] \\
&= & [(0,0)]\in\frac{\mathbb{Z}\oplus\mathbb{Z}}{\langle(2,0),(0,2)\rangle}.\\
\end{array}$$
Note that, if $f^{n}$ can be deformed to a fixed point free map over $S^{1}$ but $f^{k}$ does not then we have $b_{3}b_{4}c_{2}\equiv_{2}1$ and $k\equiv_{4}0$.

\bigskip

Case V. If $n$ is an even positive integer and $f^{n}: MA \to MA$ is deformed to a fixed point free map over $S^{1}$, then 

\noindent if $a_{3}$ is odd then $\frac{b_{4}-1}{2}((1+c_{2})(1+(n-1)b_{4})) \equiv_{2} 0 $ and at least one of $(1+c_{2})(1+(n-1)b_{4})$ and $\frac{b_{4}-1}{L}$ is even, where $L:= mdc(b_{4}-1,c_{2})$, or

\noindent if $a_{3}$ is even then $n(b_{4}-1)(c_{1}-\frac{a_{3}}{2}c_{2}-1) + (n-1)(b_{4}-1)\equiv_{2} 0$ and at least one of $n(c_{1}-\frac{a_{3}}{2}c_{2}-1) + (n-1)$ and $\frac{b_{4}-1}{L}$ is even, where $L:=mdc(b_{4}-1,c_{2}).$ 

Let $a_{3}$ odd  and $k$ an even positive integer then
$$\begin{array}{crcl}
& (1+(k-1)b_{4})) & \equiv_{2} & (1+(n-1)b_{4})) \\
\Rightarrow & \frac{b_{4}-1}{2}((1+c_{2})(1+(k-1)b_{4})) & \equiv_{2} & \frac{b_{4}-1}{2}((1+c_{2})(1+(n-1)b_{4})) \\
& & \equiv _{2}& 0 ; \\
& (1+c_{2})(1+(k-1)b_{4}) & \equiv_{2} & (1+c_{2})(1+(n-1)b_{4}).
\end{array}$$
Then, $f^{k}: MA \to MA$ can be deformed to a fixed point free map over $S^{1}$ for $a_{3}$ odd. Let $a_{3}$ even  and $k$ an even positive integer then
$$\begin{array}{crcl}
 & n(b_{4}-1)(c_{1}-\frac{a_{3}}{2}c_{2}-1) + (n-1)(b_{4}-1) & \equiv_{2} &b_{4}-1;\\
& n(c_{1}-\frac{a_{3}}{2}c_{2}-1) + (n-1) & \equiv_{2} & 1;\\
\Rightarrow & k(b_{4}-1)(c_{1}-\frac{a_{3}}{2}c_{2}-1) + (k-1)(b_{4}-1) & \equiv_{2} & 0.
\end{array}$$
Then, $f^{k}: MA \to MA$ can be deformed to a fixed point free map over $S^{1}$ for $a_{3}$ even.

\qed

\bigskip

Given $n \in \mathbb{N}$ and $f: MA \to MA$ a fiber-preserving then from Propositions \ref{proposition-3} and \ref{proposition-4} 
the conditions to deform $f$ and $f^{2}$ to a fixed point free map over $S^{1}$ are enough to deform $f^{k}$ 
to a fixed point free map over $S^{1}$ for all $k$ divisor of $n.$ 

\bigskip

\begin{theorem} \label{theorem-2} Let $f: T \times I \to T \times I$ be 
the map defined by; $$f(x,y,t)=( x+b_{3}y+c_{1}t+\varepsilon , b_{4}y+c_{2}t+\delta , t).$$ Denoting 
$f^{n}(x,y,t) = (x_{n},y_{n},t)$ then $x_{n}$ and $y_{n}$ are given by
\begin{displaymath}
	\begin{array}{ccl}
	x_{n} & = & x +b_{3}y\displaystyle{\sum_{i=0}^{n-1}b_{4}^{i}} + (nc_{1}+b_{3}c_{2}\sum_{i=0}^{n-1}ib_{4}^{n-1-i})t+b_{3}\delta\sum_{i=0}^{n-1}ib_{4}^{n-1-i}+n\varepsilon\\
	y_{n} & = & b_{4}^{n}y+c_{2}t\displaystyle{\sum_{i=0}^{n-1}b_{4}^{i}}+\delta\sum_{i=0}^{n-1}b_{4}^{i},
	\end{array}
\end{displaymath}
and for each $n \in \mathbb{N}$ and $\epsilon, \delta$ appropriates the map $f^{n}$ induces a fiber-preserving map in the fiber bundle $MA$, as in Theorem \ref{daci1theorem}, which we will represent by $f^{n}(<x,y,t>) = <x_{n},y_{n},t>,$ such that the induces homomorphism 
$(f^{n})_{\#}$ is as in the Proposition \ref{proposition-1}. Note that the induce homomorphism $f_{\#}$ is given by; 
$f_{\#}(a) = a $, $f_{\#}(b) = a^{b_{3}} b^{b_{4}} $, $f_{\#}(c) = a^{c_{1}} b^{c_{2}} c $.   \end{theorem} 

{\it Proof.} Denote $f^{n}(x,y,t)= ( x_{n},y_{n},t)$ for each $n \in \mathbb{N}.$ We have
\begin{displaymath}
	\begin{array}{ccl}
	x_{2} & = & x_{1}+b_{3}y_{1}+c_{1}t+\varepsilon\\
	      & = & (x+b_{3}y+c_{1}t+\varepsilon )+b_{3}(b_{4}y+c_{2}t+\delta)+c_{1}t+\varepsilon\\
	      & = & x +b_{3}y(b_{4}+1) + (2c_{1}+b_{3}c_{2})t+b_{3}\delta+2\varepsilon;\\
	y_{2} & = & b_{4}y_{1}+c_{2}t+\delta\\
	      & = & b_{4}(b_{4}y+c_{2}t+\delta)+c_{2}t+\delta\\     
	      & = &  b_{4}^{2}y+c_{2}(b_{4}+1)t+(b_{4}+1)\delta.
	\end{array}
\end{displaymath}
               
Suppose that $f^{n}(x,y,t) = (x_{n},y_{n},t)$ as in hypothesis, then 
$$f^{n+1}(x,y,t)=( x_{n}+b_{3}y_{n}+c_{1}t+\varepsilon , b_{4}y_{n}+c_{2}t+\delta , t)=(x_{n+1},y_{n+1},t),$$
where 
\begin{displaymath}
	\begin{array}{ccl}
	x_{n+1} & = & x_{n}+b_{3}y_{n}+c_{1}t+\varepsilon\\
	        & = & (x +b_{3}y\displaystyle{\sum_{i=0}^{n-1}b_{4}^{i}} + (nc_{1}+b_{3}c_{2}\sum_{i=0}^{n-1}ib_{4}^{n-1-i})t+
				b_{3}\delta\sum_{i=0}^{n-1}ib_{4}^{n-1-i}+n\varepsilon) + \\
          &		& + b_{3}(b_{4}^{n}y+c_{2}t\displaystyle{\sum_{i=0}^{n-1}b_{4}^{i}}+\delta\sum_{i=0}^{n-1}b_{4}^{i})+c_{1}t+\varepsilon\\
	        & = & \displaystyle{x + ( b_{3}y\sum_{i=0}^{n-1}b_{4}^{i} + b_{3}yb_{4}^{n} ) + ((nc_{1}+b_{3}c_{2}\sum_{i=0}^{n-1}ib_{4}^{n-1-i})t +c_{1}t+} \\
					&   & \displaystyle{ b_{3}c_{2}t\displaystyle{\sum_{i=0}^{n-1}b_{4}^{i}}) } \displaystyle{+ (b_{3}\delta\sum_{i=0}^{n-1}ib_{4}^{n-1-i}+ b_{3}\delta\sum_{i=0}^{n-1}b_{4}^{i}) +(n\varepsilon+\varepsilon)} \\
					& = & x +b_{3}y\displaystyle{\sum_{i=0}^{n}b_{4}^{i}} + ((n+1)c_{1}+b_{3}c_{2}\sum_{i=0}^{n}ib_{4}^{n-i})t+b_{3}\delta\sum_{i=0}^{n}ib_{4}^{n-i}+(n+1)\varepsilon;\\
\end{array}
\end{displaymath}

\begin{displaymath}
	\begin{array}{ccl}
	y_{n+1} & = & b_{4}y_{n}+c_{2}t+\delta\\
	        & = & b_{4}(b_{4}^{n}y+c_{2}t\displaystyle{\sum_{i=0}^{n-1}b_{4}^{i}}+\delta\sum_{i=0}^{n-1}b_{4}^{i})+c_{2}t+\delta\\     
	        & = & \displaystyle{b_{4}^{n+1}y+ (c_{2}t\sum_{i=1}^{n}b_{4}^{i} +c_{2}t)+(\delta\sum_{i=1}^{n}b_{4}^{i}+\delta)}\\
					& = & b_{4}^{n+1}y+c_{2}t\displaystyle{\sum_{i=0}^{n}b_{4}^{i}}+\delta\sum_{i=0}^{n}b_{4}^{i},
	\end{array}
\end{displaymath}
as we wish. 
Now, to verify that $f^{n}(<x,y,0>) = f^{n}(<A\left(^{x}_{y}\right),1>)$ for each $n \in \mathbb{N}$, in 
$T \times I$, where $T$ is the torus, firstly we will verify this condition for $n = 1$. We have
$$f<x,y,0>=< x+b_{3}y+\varepsilon , b_{4}y+\delta , 0> \,\,\,\,\, and $$
$$ f<A\left(^{x}_{y}\right),1> = <(a_{1}+a_{2}b_{3})x +(a_{3}+b_{3}a_{4})y +c_{1}+\varepsilon, 
b_{4}a_{2}x+b_{4}a_{4}y+c_{2} +\delta, 1 >$$
But in $MA$ we have $<x,y,0> = <A\left(^{x}_{y}\right),1>$, that is, 
$<x,y,0> = <a_{1}x+a_{3}y, a_{2}x+a_{4}y,1>.$
Now we will analyze each case of the Theorem \ref{daci1theorem}. 

\bigskip

Case I. \,\, In this case we need consider $b_{3} = 0$ and $b_{4}=1.$ Thus, in $MA$
$ f<x,y,0> = <x+\epsilon, y + \delta,0> = $ 
$<a_{1}x + a_{3}y + a_{1}\epsilon+a_{3}\delta, a_{2}x + a_{4}y + a_{2}\epsilon+a_{4}\delta,1>. $
Note that, $f<A\left(^{x}_{y}\right),1> = 
<a_{1}x +a_{3}y +c_{1}+\epsilon, a_{2}x+a_{4}y+c_{2} +\delta, 1 > . $
Therefore $f<x,y,0> = f<A\left(^{x}_{y}\right),1>$ if $a_{1}\epsilon+a_{3}\delta = \epsilon + k$ and 
$a_{2}\epsilon+a_{4}\delta = \delta + l$ where $k,l \in \mathbb{Z}.$

\bigskip

Case II. \,\, In this case we have $a_{1} = a_{4}=1$, $a_{2}=0$ and $a_{3}(b_{4}-1)=0.$ Therefore, 
$f<A\left(^{x}_{y}\right),1> = <x +(a_{3}+b_{3})y +c_{1}+\epsilon, b_{4}y+c_{2} +\delta, 1 >.$
Thus, $f<x,y,0> = <x+b_{3}y+\epsilon,b_{4}y+\delta,0> = <x +(a_{3}+b_{3})y +\epsilon+a_{3}\delta, b_{4}y+\delta,1>$ 
$= <x +(a_{3}+b_{3})y +\epsilon+ a_{3}\delta, b_{4}y+\delta,1>.$  Therefore 
$f<x,y,0> = f<A\left(^{x}_{y}\right),1>$ if $a_{3}\delta \in \mathbb{Z}.$

\bigskip

Case III. \,\, In this case we have $a_{1} = 1,$ $ a_{4}=-1$, $a_{2}=0$ and $a_{3}(b_{4}-1)=-2b_{3}.$ 
Therefore $f<A\left(^{x}_{y}\right),1> = <x +(a_{3}-b_{3})y +c_{1}+\epsilon, -b_{4}y+c_{2} +\delta, 1 >.$
We have $f<x,y,0> = <x+b_{3}y+\epsilon,b_{4}y+\delta,0>  = <x+(a_{3}b_{4}+b_{3})y+\epsilon +a_{3}\delta,-b_{4}y-\delta,1>.$
 Thus $f<x,y,0> = <x+(a_{3}-b_{3})y+\epsilon +a_{3}\delta,-b_{4}y-\delta,1>.$ Then 
$f<x,y,0> = f<A\left(^{x}_{y}\right),1>$ if $a_{3}\delta \in \mathbb{Z}$ and $\delta = \frac{k}{2}, k \in \mathbb{Z}.$

\bigskip

Case IV. \,\, In this case we have $a_{1} = -1,$ $ a_{4}=-1$, $a_{2}=0$ and $a_{3}(b_{4}-1)=0.$
Thus $f<A\left(^{x}_{y}\right),1> = <-x +(a_{3}-b_{3})y +c_{1}+\epsilon, -b_{4}y+c_{2} +\delta, 1 >.$
We have $f(x,y,0) = (x+b_{3}y+\epsilon,b_{4}y+\delta,0)  = (-x+(a_{3}b_{4}-b_{3})y-\epsilon +a_{3}\delta,-b_{4}y-\delta,1)$
Thus $f<x,y,0> = <-x+(a_{3}-b_{3})y-\epsilon +a_{3}\delta,-b_{4}y-\delta,1>.$ Then 
$f<x,y,0> = f<A\left(^{x}_{y}\right),1>$ if $\epsilon = \frac{a_{3}m +2k}{4}$ and $ \delta = \frac{m}{2}$ 
where $m,k \in \mathbb{Z}.$

\bigskip

Case V. \,\, In this case we have $a_{1} = -1,$ $ a_{4}=1$, $a_{2}=0$ and $a_{3}(b_{4}-1)= 2b_{3}.$
Therefore $f<A\left(^{x}_{y}\right),1> = <-x +(a_{3}+b_{3})y +c_{1}+\epsilon, b_{4}y+c_{2} +\delta, 1 >.$
We have $f(x,y,0) = (x+b_{3}y+\epsilon,b_{4}y+\delta,0)  = (-x+(a_{3}b_{4}-b_{3})y-\epsilon +a_{3}\delta,b_{4}y+\delta,1)$
 Thus $f<x,y,0> = <x+(a_{3}+b_{3})y-\epsilon +a_{3}\delta,b_{4}y+\delta,1>.$ Then 
$f<x,y,0> = f<A\left(^{x}_{y}\right),1>$ if $ \epsilon = \frac{a_{3}\delta+k}{2}$  
where $k \in \mathbb{Z}.$ 

Now we will verify the condition for all $n \in \mathbb{N}.$ 

Case I. In $MA$ we have
$$\begin{array}{ccl}	
f^{n}<x,y,0> & = & <x +n\varepsilon, y +n\delta,0>\\
             & = & <a_{1}(x+n\varepsilon) + a_{3}(y + n\delta), a_{2}(x+n\varepsilon) + a_{4}(y+n\delta) ,1>\\
				 & = & <a_{1}x + a_{3}y + na_{1}\varepsilon+na_{3}\delta, a_{2}x + a_{4}y + na_{2}\varepsilon+na_{4}\delta,1>.
\end{array}$$
But, $f^{n}\left<A \left(^{x}_{y}\right),1\right>	= <a_{1} x + a_{3} y + n c_{1} +n\varepsilon, a_{2}x + a_{4}y + nc_{2}+ n\delta,1>.$
Then, $f^{n}<x,y,0> = f^{n}\left<A \left(^{x}_{y}\right),1\right>$ if $na_{1}\varepsilon +n a_{3}\delta= n\varepsilon+k$, $k\in\mathbb{Z}$, and $na_{2}\varepsilon+na_{4}\delta = n\delta +l$, $l \in\mathbb{Z}$.

\bigskip

Case II. 
$$\begin{array}{l}	
f^{n}<x,y,0> = \\
=  <x +b_{3}y\displaystyle{\sum_{i=0}^{n-1}b_{4}^{i}} +b_{3}\delta\sum_{i=0}^{n-1}ib_{4}^{n-1-i}+n\varepsilon, b_{4}^{n}y+\delta\sum_{i=0}^{n-1}b_{4}^{i},0 >\\
				  =  <\left(x +b_{3}y\displaystyle{\sum_{i=0}^{n-1}b_{4}^{i}} +b_{3}\delta\sum_{i=0}^{n-1}ib_{4}^{n-1-i}+n\varepsilon\right) + a_{3}\left( b_{4}^{n}y+\delta\displaystyle{\sum_{i=0}^{n-1}b_{4}^{i}}\right) , \\ 
					b_{4}^{n}y+ \delta\displaystyle{\sum_{i=0}^{n-1}b_{4}^{i}},1>\\
				   =  <x + \left(a_{3}b_{4}^{n}+b_{3}\displaystyle{\sum_{i=0}^{n-1}b_{4}^{i}}     \right)y + 
					\left(a_{3}\displaystyle{\sum_{i=0}^{n-1}b_{4}^{i}}+b_{3}\displaystyle{\sum_{i=0}^{n-1}ib_{4}^{n-1-i}}     \right)\delta + 
					n\varepsilon, \\
					b_{4}^{n}y+\delta\displaystyle{\sum_{i=0}^{n-1}b_{4}^{i}},1 >
\end{array}$$
But, $f^{n}\left(<A \left(^{x}_{y}\right),1>\right)=$
$$\begin{array}{l}
 	=<x + \left(a_{3}+b_{3}\displaystyle{\sum_{i=0}^{n-1}b_{4}^{i}}\right)y + 
	b_{3}\delta\displaystyle{\sum_{i=0}^{n-1}ib_{4}^{n-1-i}}+ n\varepsilon+ 
	b_{3}c_{2}\displaystyle{\sum_{i=0}^{n-1}ib_{4}^{n-1-i}}  
	+ nc_{1}, \\ b_{4}^{n}y+\delta\displaystyle{\sum_{i=0}^{n-1}b_{4}^{i}} + c_{2}\displaystyle{\sum_{i=0}^{n-1}b_{4}^{i}},1>.
\end{array}$$
Thus $f^{n}\left(<A \left(^{x}_{y}\right),1>\right) = f^{n}<x,y,0> $ if 
$\delta a_{3} \displaystyle{\sum_{i=0}^{n-1} b_{4}^{i}}\in\mathbb{Z}$. 

\bigskip

Case III. $f^{n}<x,y,0>  = $
$$\begin{array}{l}	
 =< x +b_{3}y\displaystyle{\sum_{i=0}^{n-1}b_{4}^{i}} + b_{3}\delta\sum_{i=0}^{n-1}ib_{4}^{n-1-i}+n\varepsilon , b_{4}^{n}y+\delta\sum_{i=0}^{n-1}b_{4}^{i},0 >\\
         =  <\left(x +b_{3}y\displaystyle{\sum_{i=0}^{n-1}b_{4}^{i}} +b_{3}\delta\sum_{i=0}^{n-1}ib_{4}^{n-1-i}+n\varepsilon\right) + a_{3}\left( b_{4}^{n}y+\delta\displaystyle{\sum_{i=0}^{n-1}b_{4}^{i}}\right) ,\\
				-b_{4}^{n}y-\delta\displaystyle{\sum_{i=0}^{n-1}b_{4}^{i}},1 >\\
				=  <x + \left(a_{3}b_{4}^{n}+b_{3}\displaystyle{\sum_{i=0}^{n-1}b_{4}^{i}}     \right)y +
				\left(a_{3}\displaystyle{\sum_{i=0}^{n-1}b_{4}^{i}}+b_{3}\displaystyle{\sum_{i=0}^{n-1}ib_{4}^{n-1-i}}   
					\right)\delta + n\varepsilon, \\
					-b_{4}^{n}y-\delta\displaystyle{\sum_{i=0}^{n-1}b_{4}^{i}},1>
\end{array}$$
 Note that,  $f^{n}\left(<A \left(^{x}_{y}\right),1>\right)	 = $
$$\begin{array}{l}
=<x + \left(a_{3}-b_{3}\displaystyle{\sum_{i=0}^{n-1}b_{4}^{i}}\right)y - b_{3}\delta\displaystyle{\sum_{i=0}^{n-1}ib_{4}^{n-1-i}}+ n\varepsilon-b_{3}c_{2}\displaystyle{\sum_{i=0}^{n-1}ib_{4}^{n-1-i}}    + nc_{1}, \\
 -b_{4}^{n}y+\delta\displaystyle{\sum_{i=0}^{n-1}b_{4}^{i}} + c_{2}\displaystyle{\sum_{i=0}^{n-1}b_{4}^{i}},1>.
\end{array}$$
Therefore, $f^{n}<x,y,0> = f^{n}\left(<A \left(^{x}_{y}\right),1>\right)	$ if $2\delta\displaystyle{\sum_{i=0}^{n-1}b_{4}^{i}} 
\in\mathbb{Z}$ and \break $\left(a_{3}\displaystyle{\sum_{i=0}^{n-1}b_{4}^{i}}+2b_{3}\displaystyle{\sum_{i=0}^{n-1}ib_{4}^{n-1-i}}
 \right)\delta \in \mathbb{Z}$.

\bigskip

Case IV. $f^{n}<x,y,0>  = $
$$\begin{array}{l}	
 = < x +b_{3}y\displaystyle{\sum_{i=0}^{n-1}b_{4}^{i}} + b_{3}\delta\sum_{i=0}^{n-1}ib_{4}^{n-1-i}+n\varepsilon , b_{4}^{n}y+\delta\sum_{i=0}^{n-1}b_{4}^{i},0>\\
         = < -\left(x +b_{3}y\displaystyle{\sum_{i=0}^{n-1}b_{4}^{i}} +b_{3}\delta\sum_{i=0}^{n-1}ib_{4}^{n-1-i}+n\varepsilon\right) + a_{3}\left( b_{4}^{n}y+\delta\displaystyle{\sum_{i=0}^{n-1}b_{4}^{i}}\right) , \\ 
				-b_{4}^{n}y-\delta\displaystyle{\sum_{i=0}^{n-1}b_{4}^{i}},1>\\
	\end{array}$$ 
	$$\begin{array}{l}			
				=  <-x + \left(a_{3}b_{4}^{n}-b_{3}\displaystyle{\sum_{i=0}^{n-1}b_{4}^{i}}     \right)y +
				\left(a_{3}\displaystyle{\sum_{i=0}^{n-1}b_{4}^{i}}-b_{3}\displaystyle{\sum_{i=0}^{n-1}ib_{4}^{n-1-i}}     \right)\delta 
				- n\varepsilon, \\ 
				-b_{4}^{n}y-\delta\displaystyle{\sum_{i=0}^{n-1}b_{4}^{i}},1>
\end{array}$$
We have,  $f^{n}\left(<A \left(^{x}_{y}\right),1>\right)	 = $
$$\begin{array}{l}						
= <-x + \left(a_{3}-b_{3}\displaystyle{\sum_{i=0}^{n-1}b_{4}^{i}}\right)y - b_{3}\delta\displaystyle{\sum_{i=0}^{n-1}ib_{4}^{n-1-i}}+ 
n\varepsilon-b_{3}c_{2}\displaystyle{\sum_{i=0}^{n-1}ib_{4}^{n-1-i}}    + nc_{1},\\
 -b_{4}^{n}y+\delta\displaystyle{\sum_{i=0}^{n-1}b_{4}^{i}} + c_{2}\displaystyle{\sum_{i=0}^{n-1}b_{4}^{i}},1 >
\end{array}$$
Thus, $f^{n}<x,y,0> = f^{n}\left(<A \left(^{x}_{y}\right),1>\right)	$ if  $2\delta\displaystyle{\sum_{i=0}^{n-1}b_{4}^{i}}\in\mathbb{Z}$ 
and $2n\varepsilon=a_{3}\delta\displaystyle{\sum_{i=0}^{n-1}b_{4}^{i}}+ k$, $k\in\mathbb{Z}$. 

\bigskip

Case V. $f^{n}<x,y,0>  = $
$$\begin{array}{l}	
= < x +b_{3}y\displaystyle{\sum_{i=0}^{n-1}b_{4}^{i}} + b_{3}\delta\sum_{i=0}^{n-1}ib_{4}^{n-1-i}+n\varepsilon , 
b_{4}^{n}y+\delta\sum_{i=0}^{n-1}b_{4}^{i},0>\\
  =  <-\left(x +b_{3}y\displaystyle{\sum_{i=0}^{n-1}b_{4}^{i}} +b_{3}\delta\sum_{i=0}^{n-1}ib_{4}^{n-1-i}+n\varepsilon\right) 
	+ a_{3}\left( b_{4}^{n}y+\delta\displaystyle{\sum_{i=0}^{n-1}b_{4}^{i}}\right) , \\ 
	b_{4}^{n}y+ 	\delta\displaystyle{\sum_{i=0}^{n-1}b_{4}^{i}},1>\\
	=  < -x + \left(a_{3}b_{4}^{n}-b_{3}\displaystyle{\sum_{i=0}^{n-1}b_{4}^{i}}     \right)y + 
	\left(a_{3}\displaystyle{\sum_{i=0}^{n-1}b_{4}^{i}}-b_{3}\displaystyle{\sum_{i=0}^{n-1}ib_{4}^{n-1-i}}  \right)\delta 
	- n\varepsilon, \\ 
	b_{4}^{n}y+\delta\displaystyle{\sum_{i=0}^{n-1}b_{4}^{i}},1 >
\end{array}$$
We have,  $f^{n}\left(<A \left(^{x}_{y}\right),1>\right)	 = $
$$\begin{array}{l}
 =<-x + \left(a_{3}+b_{3}\displaystyle{\sum_{i=0}^{n-1}b_{4}^{i}}\right)y - b_{3}\delta\displaystyle{\sum_{i=0}^{n-1}ib_{4}^{n-1-i}}+ n\varepsilon-b_{3}c_{2}\displaystyle{\sum_{i=0}^{n-1}ib_{4}^{n-1-i}}    + 
nc_{1}, \\ 
+b_{4}^{n}y+\delta\displaystyle{\sum_{i=0}^{n-1}b_{4}^{i}} + c_{2}\displaystyle{\sum_{i=0}^{n-1}b_{4}^{i}},1>.
\end{array}$$
Therefore, $f^{n}<x,y,0> = f^{n}\left(<A \left(^{x}_{y}\right),1>\right)	$ if $2n\varepsilon=a_{3}\delta\displaystyle{\sum_{i=0}^{n-1}b_{4}^{i}}+ k$, $k\in\mathbb{Z}$.

Thus for each $n \in \mathbb{N}$ and $\epsilon, \delta $ satisfying the conditions above the map $f^{n}: T \times I \to T \times I$ 
induces a fiber-preserving map on $MA$ which will be represent by the same symbol. \qed

\bigskip

%In Theorem \ref{theorem-2} in all cases for $\epsilon = \delta = 0$ the map $f$ induces a fiber-preserving map on $MA,$ 
%that is, $f<x,y,t>=< x+b_{3}y+c_{1}t+ , b_{4}y+c_{2}t , t>$ is a fiber-preserving map on $MA.$ Now, given $\epsilon, \delta \in \mathbb{R}$  %consider the homotopy $H: MA \times I \to MA$ defined by; 
%$$H(<x,y,t>,s) = < x+b_{3}y+c_{1}t+ s\epsilon , b_{4}y+c_{2}t + s\delta , t> $$
%Note that $p \circ H = p.$ Therefore each fiber-preserving map $f<x,y,t>=< x+b_{3}y+c_{1}t , b_{4}y+c_{2}t , t>$ is fiberwise homotopy 
%to the map $\overline{f}<x,y,t> = < x+b_{3}y+c_{1}t+ \epsilon , b_{4}y+c_{2}t + \delta , t>$ for any $\epsilon, \delta \in \mathbb{R}.$
%In this sense, sometimes we use a fiber-preserving map $f^{n}$ as in the Theorem \ref{theorem-2} or a map fiberwise homotopic to it. 
%
The next result we will help us to study the fixed points of $f^{n}$ given in Theorem \ref{theorem-2}. 

\begin{proposition}\label{prop-2}
Let $n, \ b_{3}, \ b_{4}, \ c_{1}, \ c_{2} \in \mathbb{Z}$, $n\geq1$. If $c_{1}(b_{4}-1)-c_{2}b_{3}\neq 0$ then for all $\varepsilon , \ \delta \in \mathbb{R}$ there are $a,\ b \in \mathbb{Z}$ such that the system bellow has solution $(x,y,t) \in \mathbb{R}^{2}\times I$:
$$ \left\{
\begin{array}{lll}
 x + a & = & x +b_{3}y\displaystyle{\sum_{i=0}^{n-1}b_{4}^{i}} + (nc_{1}+b_{3}c_{2}\sum_{i=0}^{n-1}ib_{4}^{n-1-i})t+b_{3}\delta\sum_{i=0}^{n-1}ib_{4}^{n-1-i}+n\varepsilon;\\
 y + b & = & b_{4}^{n}y+c_{2}t\displaystyle{\sum_{i=0}^{n-1}b_{4}^{i}}+\delta\sum_{i=0}^{n-1}b_{4}^{i}. 
\end{array} \right.
$$ 
\end{proposition}

{\it Proof:} Suppose $b_{4}\neq 1$ and $b_{4}\neq-1$ with $n$ even ($b_{4}=-1$ with $n$ odd is allowed) and $c_{1}(b_{4}-1)-b_{3}c_{2}\neq 0$ then given $\varepsilon, \ \delta\in\mathbb{R}$ we have the solutions $x\in\mathbb{R}$ and:
\begin{displaymath}
	\begin{array}{ccl}
		t & = & \dfrac{nb_{3}\delta-n(b_{4}-1)\varepsilon - (b_{4}-1)a - b_{3}b}{n(c_{1}(b_{4}-1)-b_{3}c_{2})}\in I;\\
		y & = & \dfrac{nc_{2}\varepsilon-nc_{1}\delta - ac_{2}}{n(c_{1}(b_{4}-1)-b_{3}c_{2})} + b\left(\dfrac{1}{b^{n}_{4}-1}+ \dfrac{b_{3}c_{2}}{n(b_{4}-1)(c_{1}(b_{4}-1)-b_{3}c_{2})}\right)\in\mathbb{R}.
	\end{array}
\end{displaymath}
Thus, we need to find $a, \ b\in \mathbb{Z}$ such that $0\leq t\leq 1$. Let $k_{0}=n(c_{1}(b_{4}-1)-b_{3}c_{2})\in \mathbb{Z}$, $k_{0}\neq0$, and $k_{1}= nb_{3}\delta-n(b_{4}-1)\varepsilon\in \mathbb{R}$, $t=\dfrac{k_{1} - (b_{4}-1)a - b_{3}b}{k_{0}}$. If $0\leq k_{1}\leq k_{0}$ or $k_{0}\leq k_{1}\leq 0$ let $a=b=0$, then $t=\dfrac{k_{1}}{k_{0}}$. If $0< k_{0}\leq k_{1}$ or $k_{1}\leq0< k_{0}$ then there are $d, \ q \in \mathbb{Z}$ such that $k_{1}=dk_{0}+q$ with $0\leq q <k_{0}$. Let $a=nc_{1}d$ and $b=nc_{2}d$, then 
$$t=\dfrac{dk_{0}+q - (b_{4}-1)nc_{1}d - b_{3}nc_{2}d}{k_{0}}=d +\dfrac{q}{k_{0}}- \dfrac{dk_{0}}{k_{0}}=\dfrac{q}{k_{0}}.$$
If $ k_{1}\leq k_{0}<0$ or $k_{0}<0\leq k_{1}$ then there are $d, \ q \in \mathbb{Z}$ such that $k_{1}=dk_{0}+q$ with $0\leq q <|k_{0}|$. Let $k\in\mathbb{Z}$ the least integer greater than $\dfrac{-q}{k_{0}}$,  $a=nc_{1}(d-k)$ and $b=nc_{2}(d-k)$, then 
$$t=\dfrac{dk_{0}+q - (b_{4}-1)nc_{1}(d-k) - b_{3}nc_{2}(d-k)}{k_{0}}=d +\dfrac{q}{k_{0}}- (d-k)=\dfrac{q}{k_{0}}+k.$$
Then, $0\leq t\leq 1$.

%

%*************************

%

If $b_{4}= 1$ and $c_{1}(b_{4}-1)-b_{3}c_{2}\neq 0$  then $b_{3}c_{2}\neq 0$. Thus, given $\varepsilon, \ \delta\in\mathbb{R}$ we have the solutions $x\in\mathbb{R}$ and:
\begin{displaymath}
	\begin{array}{ccl}
		t & = & \dfrac{b}{nc_{2}}-\dfrac{\delta}{c_{2}}\in I;\\
		y & = & \dfrac{-nc_{2}\varepsilon+nc_{1}\delta + ac_{2}}{nb_{3}c_{2}} - b\left(\dfrac{c_{1}}{nb_{3}c_{2}}+ \dfrac{n-1}{2n}\right)\in\mathbb{R}.
	\end{array}
\end{displaymath}
We need to find $b\in \mathbb{Z}$ such that $0\leq t\leq 1$. If $c_{2}>0$ take $n\delta\leq b\leq n(c_{2}+\delta)$ and if $c_{2}<0$ take $n\delta\geq b\geq n(c_{2}+\delta)$.  

%

%*************************

%

If $b_{4}= -1$, $n$ even and $c_{1}(b_{4}-1)-b_{3}c_{2}\neq 0$ then $2c_{1}+ b_{3}c_{2}\neq 0$. Thus, given $\varepsilon, \ \delta\in\mathbb{R}$ we have the solutions $x,\ y\in\mathbb{R}$ and:
\begin{displaymath}
	\begin{array}{ccl}
		t & = & \dfrac{2a-nb_{3}\delta-2n\varepsilon}{n(2c_{1}+ b_{3}c_{2})}\in I;\\
	\end{array}
\end{displaymath}
We need to find $a\in \mathbb{Z}$ such that $0\leq t\leq 1$. Let $n=2k$, $k_{0}=k(2c_{1}+b_{3}c_{2})\in \mathbb{Z}$, $k_{0}\neq0$, and $k_{1}= -kb_{3}\delta-2k\varepsilon\in \mathbb{R}$, $t=\dfrac{a+k_{1}}{k_{0}}$. If $0\leq k_{1}\leq k_{0}$ or $k_{0}\leq k_{1}\leq 0$ let $a=0$, then $t=\dfrac{k_{1}}{k_{0}}$. If $0< k_{0}\leq k_{1}$ or $k_{1}\leq0< k_{0}$ then there are $d, \ q \in \mathbb{Z}$ such that $k_{1}=dk_{0}+q$ with $0\leq q <k_{0}$. Let $a=-k_{0}d$, then $t=\dfrac{-k_{0}d + dk_{0}+q}{k_{0}}=\dfrac{q}{k_{0}}$. If $ k_{1}\leq k_{0}<0$ or $k_{0}<0\leq k_{1}$ then there are $d, \ q \in \mathbb{Z}$ such that $k_{1}=dk_{0}+q$ with $0\leq q <|k_{0}|$. Let $l\in\mathbb{Z}$ the greatest integer lower than $q$ and  $a=-l-k_{0}d$, then $t=\dfrac{-l-k_{0}d+dk_{0}+q}{k_{0}}=\dfrac{q-l}{k_{0}}$. Then, $0\leq t\leq 1$.

\qed

%*******************************************************

\begin{theorem}[Main Theorem] \label{maintheorem}
Let $f: MA \to MA$ be a fiber-preserving map, where $MA$ is a T-bundle over $S^{1},$ 
and $f_{\#}(a) = a $, $f_{\#}(b) = a^{b_{3}} b^{b_{4}} $, $f_{\#}(c) = a^{c_{1}} b^{c_{2}} c $. 
Suppose that $f^{n}: MA \to MA$ can be deformed to a fixed point free map over $S^{1}.$
Then there exits a map $g$ fiberwise homotopic to $f$ such that $g^{n}$ is a fixed point free map, in the 
cases $I$ and $II,$ if and only if the following conditions are satisfies;

\bigskip
1) $MA$ is as in case $I$ and $f$ is arbitrary.

2) $MA$ is as in cases $II$, $III$ and  $c_{1}(b_{4}-1)-c_{2}b_{3} = 0$.

\end{theorem}
{\it Proof.}
(1) For each map $f$ such that $(f_{|T})_{\#} = Id$ consider the map $g$ fiberwise homotopic to $f$ given by;
$g^{'}(<x,y,t>) = <x+c_{1}t+\epsilon, y+c_{2}t+\delta, t>$, with $\epsilon , \delta \in \mathbb{Q}-\mathbb{Z} $ 
satisfying the condition; $a_{1}\epsilon + a_{3}\delta = \epsilon + k$ and 
$a_{2}\epsilon + a_{4}\delta = \epsilon + k$ for some $k,l \in \mathbb{Z}.$  Note that 
$g^{'}$ is fiberwise homotopic to the map $g$ defined by;
$$ g(<x,y,t>) = \left \{
\begin{array}{lll}
<x+2c_{1}t+\epsilon, y + \delta, t> & if & 0 \leq t \leq \frac{1}{2} \\
<x+ c_{1}+\epsilon, y + c_{2}(2t-1)+\delta, t> & if & \frac{1}{2} \leq t \leq 1 \\
\end{array} \right.
$$
In fact, $H: MA \times I \to MA$ defined by; 
$$ H(<x,y,t>,s) = \left \{
\begin{array}{lll}
<x+c_{1}t+\epsilon, y + c_{2}t + \delta, t> & if & 0 \leq t \leq s \\
<x+ c_{1}(2t-s)+\epsilon, y + c_{2}s+\delta, t> & if & s \leq t \leq \frac{s+1}{2} \\
<x+ c_{1}+\epsilon, y + c_{2}(2t-1)+\delta, t> & if & \frac{s+1}{2} \leq t \leq 1 \\
\end{array} \right.
$$
is a homotopy between $g^{'}$ and $g.$ Note that, 
$$ g^{n}(<x,y,t>) = \left \{
\begin{array}{lll}
<x+n2c_{1}t+n\epsilon, y + n\delta, t> & if & 0 \leq t \leq \frac{1}{2} \\
<x+ nc_{1}+n\epsilon, y + nc_{2}(2t-1)+n\delta, t> & if & \frac{1}{2} \leq t \leq 1 \\
\end{array} \right.
$$ 
In this case we have $Fix(g^{n}) = \emptyset$. 

\bigskip

(2) \,\, From Theorem \ref{theorem-1} we must consider two situations; $c_{1}(b_{4}-1)-c_{2}b_{3} = 0$ or 
$n$ even and $b_{4} = -1.$ First we will suppose $c_{1}(b_{4}-1)-c_{2}b_{3} = 0$. Note that if $c_{1}(b_{4}-1)-c_{2}b_{3} \neq 0$
then $g$ can not be deformed to a fixed point free map. Therefore can not exits $f$ fiberwise homotopic to $g$ such that 
$Fix(f^{n}) = \emptyset$ for all $n \geq 1$ by Lemma \ref{lemma1}.

\bigskip
Suppose $b_{4}=1$ and $f^{n}(x,y,t) = (x_{n},y_{n},t)$ for each $n \in \mathbb{N}$ gives in the Theorem \ref{theorem-2}. Then:
\begin{displaymath}
	\left\{\begin{array}{ccl}
	x_{n} & = & \displaystyle{x +nb_{3}y + \left(c_{1}+b_{3}c_{2}\frac{n-1}{2}\right)nt+\frac{n(n-1)}{2}b_{3}\delta+n\varepsilon},\\
	y_{n} & = & \displaystyle{y+nc_{2}t+n\delta}.
	\end{array}\right.
\end{displaymath}
For $c_{2}b_{3} = 0$, $f^{n}(x,y,t)$ is a fixed point free map for each $n$ choosing $\epsilon \in \mathbb{R}-\mathbb{Q}$ and $\delta = 0$ if $b_{3} = 0$ or $\delta \in \mathbb{R}-\mathbb{Q}$ if $c_{2} = 0.$ 
 In fact, if $b_{3} = 0$ and $c_{2} \neq 0$ choose $\epsilon \in \mathbb{R}-\mathbb{Q}$ and $\delta = 0$, then 
$$ \left\{
\begin{array}{lll}
 x + k_{n} & = & x + nc_{1}t+ n\varepsilon  \\
 y + l_{n} & = & y + nc_{2}t 
\end{array} \right.
$$ 
for some $k_{n},l_{n} \in \mathbb{Z}.$ From equations we obtain $t = \dfrac{l_{n}}{nc_{2}} \in \mathbb{Q}$ and $ \epsilon = \dfrac{k_{n}}{n} - \dfrac{c_{1}l_{n}}{nc_{2}} \in \mathbb{Q}$, but this is a contradiction because $\epsilon $ in $\mathbb{R}-\mathbb{Q}.$ Therefore $f^{n}$ has no fixed point for all $n$. 
If $c_{2} = 0$ we choose $\delta \in \mathbb{R}-\mathbb{Q} $ then $y +n\delta = y+l_{n}$ which implies $\delta = \frac{l_{n}}{n} \in \mathbb{Q}$ that is a 
contradiction because $\delta $  in $ \mathbb{R}-\mathbb{Q} $. Thus $f^{n}$ has no fixed point for all $n$. 

\bigskip

Now we suppose $b_{4}\neq1$, $\delta = 0$ and $c_{1}(b_{4} - 1)=c_{2}b_{3}$. Then  
\begin{displaymath}
	\begin{array}{ccl}
	x_{n} & = & x +b_{3}y\displaystyle{\sum_{i=0}^{n-1}b_{4}^{i}} + (nc_{1}+b_{3}c_{2}\sum_{i=0}^{n-1}ib_{4}^{n-1-i})t+n\varepsilon\\
	      & = & x + \left( \dfrac{b_{3}(b_{4}^{n}-1)}{b_{4}-1}\right) y + \left( n c_{1}+\dfrac{b_{3}c_{2}(b_{4}^{n}-1+n(1-b_{4}))}{(b_{4}-1)^{2}}\right) t + n\varepsilon\\
				& = & x + \left( \dfrac{b_{4}^{n}-1}{b_{4}-1}\right)b_{3} y + \left( \dfrac{b_{4}^{n}-1}{b_{4}-1}\right) c_{1} t + n\varepsilon;\\
	y_{n} & = & b_{4}^{n}y+c_{2}t\displaystyle{\sum_{i=0}^{n-1}b_{4}^{i}}\\
	      & = & b_{4}^{n}y+ \left( \dfrac{c_{2}(b_{4}^{n}-1)}{b_{4}-1}\right)t.
	\end{array}
\end{displaymath}

We are interested at the solutions of the system above: 
\begin{eqnarray}
x_{n} & = & x + k_{n}; \label{eq1}\\
y_{n} & = & y + l_{n}, \label{eq2}
\end{eqnarray}  
where $k_{n}, l_{n}\in \mathbb{Z}$ and $n\in\mathbb{N}$.
With the equation \ref{eq2} and $c_{2}\neq 0$, we have:
$$ t=\dfrac{l_{n}(b_{4}-1)}{c_{2}(b_{4}^{n}-1)} + \dfrac{1-b_{4}}{c_{2}}y.$$ 
For this $t$ we obtain:
\begin{displaymath}
	\begin{array}{ccl}
	x_{n} & = &  x + \left( \dfrac{b_{4}^{n}-1}{b_{4}-1}\right)b_{3} y + \left( \dfrac{b_{4}^{n}-1}{b_{4}-1}\right) c_{1} t + n\varepsilon\\	
        & = &  x -\left( \dfrac{(b_{4}^{n}-1) ((b_{4}-1)c_{1}-c_{2}b_{3})}{(b_{4}-1) c_{2}}\right) y + n\varepsilon +\dfrac{c_{1}l_{n}}{c_{2}}\\
				& = &  x + n\varepsilon +\dfrac{c_{1}l_{n}}{c_{2}}
	\end{array}
\end{displaymath}

But, for $\varepsilon\in\mathbb{R}- \mathbb{Q}$ the equation \ref{eq2} has no solution for all $n\in\mathbb{N}$. In fact, for any $k_{n}\in\mathbb{Z}$ we have:
$$x_{n} = x + k_{n} 
\Rightarrow   x + n\varepsilon +\dfrac{c_{1}l_{n}}{c_{2}} = x + k_{n} 
\Rightarrow \underbrace{\varepsilon}_{\in\mathbb{R}-\mathbb{Q}}= \underbrace{\dfrac{k_{n}}{n} - \dfrac{c_{1}l_{n}}{nc_{2}}}_{\in\mathbb{Q}}.$$

In the other side, if $c_{2}= 0$ then $c_{1}=0$ because $b_{4}\neq 1$ and: 
\begin{displaymath}
	\begin{array}{ccl}
	x_{n} & = &  x + \left( \dfrac{b_{4}^{n}-1}{b_{4}-1}\right)b_{3} y + \left( \dfrac{b_{4}^{n}-1}{b_{4}-1}\right) c_{1} t + n\varepsilon\\
	      & = &  x + \left( \dfrac{b_{4}^{n}-1}{b_{4}-1}\right)b_{3} y +  n\varepsilon;\\
	y_{n} & = & b_{4}^{n}y+ \left( \dfrac{c_{2}(b_{4}^{n}-1)}{b_{4}-1}\right)t\\
	      & = & b_{4}^{n}y.
	\end{array}
\end{displaymath}

Thus, $y=\dfrac{l_{n}}{b_{4}^{n}-1}$ by equation \ref{eq2}. Then, $x_{n} = x + n\varepsilon + \dfrac{b_{3}l_{n}}{b_{4}-1}$ and the equation \ref{eq1} has no solution for $\varepsilon\in\mathbb{R}- \mathbb{Q}$ and for $n\in\mathbb{N}$. In fact, for any $k_{n}, l_{n}\in\mathbb{Z}$ we have:
$$x_{n} = x + k_{n} 
\Rightarrow   x + n\varepsilon +\dfrac{b_{3}l_{n}}{b_{4}-1} = x + k_{n} 
\Rightarrow \underbrace{\varepsilon}_{\in\mathbb{R}-\mathbb{Q}}= \underbrace{\dfrac{k_{n}}{n} - \dfrac{b_{3}l_{n}}{n(b_{4}-1)}}_{\in\mathbb{Q}}.$$

If $b_{4}=- 1$ and $n$ even then $f^{n}(x,y,t) = (x_{n},y_{n},t)$ such that: 
\begin{displaymath}
	\begin{array}{ccl}
	x_{n} & = & x + b_{3}\delta\left(\frac{n}{2}\right)+n\varepsilon;\\
	y_{n} & = & y.
	\end{array}
\end{displaymath}
Thus, choosing $\delta = 0$ and $\varepsilon \in \mathbb{R}-\mathbb{Q},$ the first equation has no solution, that is, $Fix(f) =  \emptyset.$
\qed

\bigskip
With a similar way of the Theorem \ref{maintheorem} we obtain the following result.

\begin{theorem}
Let $f: MA \to MA$ be a fiber-preserving map, where $MA$ is a T-bundle over $S^{1},$ 
and $f_{\#}(a) = a $, $f_{\#}(b) = a^{b_{3}} b^{b_{4}} $, $f_{\#}(c) = a^{c_{1}} b^{c_{2}} c $. 
Suppose that $f^{n}: MA \to MA$ can be deformed to a fixed point free map over $S^{1}$ as in the Theorem \ref{theorem-1}. 
If the following conditions are satisfied, in the cases $IV$ and $V,$ then there exits a map $g$ fiberwise homotopic to $f$ such that 
$g^{n}$ is a fixed point free map.

\bigskip

$MA$ is as in case $IV,$ $V$ and 
$$c_{1}(b_{4}-1)-c_{2}b_{3} = 0 \,\, and \,\, n \,\, odd \,\, or $$ 
$$c_{1}(b_{4}-1)-c_{2}b_{3} = 0, \,\, b_{4} \,\, odd \,\, and \,\, n = 4k + 2, k \geq 0. $$

\end{theorem}

%{\it Proof.}
%(1) Let $n$ be a positive integer and  $g: MA \to MA$ be 
%the map defined by:
% $$g(<x,y,t>)=< x+b_{3}y+c_{1}t+\varepsilon , b_{4}y+c_{2}t+\delta , t>,$$
%such that $c_{1}(b_{4}-1)-c_{2}b_{3} = 0.$ Then for $b_{3}\neq0$ and a positive divisor $k$ of $n$ the system of Proposition \ref{prop-2} for %$g^{k}$ is such that $b  =   \dfrac{(b_{4}-1)(a-k \varepsilon)}{b_{3}} + k\delta$. Note that $b\notin\mathbb{Z}$ for $\varepsilon=0$ and %$\delta\in\mathbb{R-Q}$ for all $k$ divisor of $n$. So, $g$ is fiberwise homotopic to $f(x,y,t)=( x+b_{3}y+c_{1}t, b_{4}y+c_{2}t+\pi , t)$ %and $f^{n}$ is a fixed point free map.
%
%If $b_{3}=0$ and $c_{1}\neq0$ then $b_{4}=1$. Thus, the system of Proposition \ref{prop-2} for $g^{k}$ is such that $b = %\dfrac{c_{2}(a-k\varepsilon)}{c_{1}}+ k\delta$. Note that $b\notin\mathbb{Z}$ when $\varepsilon=0$ and $\delta\in\mathbb{R-Q}$ for all $k$ %divisor of $n$.
%Therefore $g$ is fiberwise homotopic to $f(x,y,t)=( x+c_{1}t, y+c_{2}t+\pi , t)$ and $f^{n}$ is a fixed point free map.
%
%Finally, if $b_{3}=c_{1}=0$ then the system of Proposition \ref{prop-2} for $g^{k}$ is such that $a  =   k\varepsilon$. In this case, %$a\notin\mathbb{Z}$ when $\varepsilon\in\mathbb{R-Q}$ for all $k$ divisor of $n$. Therefore $g$ is fiberwise homotopic to $f(x,y,t)=( x+\pi, %b_{4}y+c_{2}t, t)$ and $f^{n}$ is a fixed point free map.
%\qed

\end{document}